\documentclass[11pt]{article}

\usepackage[T1]{fontenc}
\usepackage[utf8]{inputenc}
\usepackage{lmodern}
\usepackage[a4paper,margin=1.1in]{geometry}
\usepackage{graphicx}
\usepackage{amsmath,amssymb,amsfonts,amsthm}
\usepackage{mathtools}
\usepackage{mathrsfs}
\usepackage{booktabs}
\usepackage{multirow}
\usepackage{makecell,array}
\usepackage{enumitem}
\usepackage{siunitx}
\usepackage{algorithm}
\usepackage{algorithmicx}
\usepackage{algpseudocode}
\usepackage{listings}
\usepackage{xcolor}
\usepackage{textcomp}
\usepackage{hyperref}

\hypersetup{
  colorlinks=false,
  linkcolor=black,
  citecolor=black,
  urlcolor=black,
  filecolor=black,
  pdfborder={0 0 0}
}

\theoremstyle{plain}
\newtheorem{theorem}{Theorem}
\newtheorem{proposition}[theorem]{Proposition}
\newtheorem{corollary}[theorem]{Corollary}
\newtheorem{lemma}[theorem]{Lemma}

\theoremstyle{definition}

\theoremstyle{remark}

\newcounter{subfig}

\title{One-dimensional optimisation of indefinite-weight principal eigenvalues with asymmetric Robin parameters and a Schrödinger-type perturbation}

\author{
Baruch Schneider\\
\small Department of Mathematics, University of Ostrava, Ostrava, Czech Republic\\
\small \texttt{baruch.schneider@osu.cz}
\and
Diana Schneiderov\'a\\
\small Department of Mathematics, University of Ostrava, Ostrava, Czech Republic\\
\small \texttt{diana.schneiderova@osu.cz}
\and
Yifan Zhang\thanks{Corresponding author.}\\
\small Department of Mathematics, University of Ostrava, Ostrava, Czech Republic\\
\small Department of Applied Mathematics, VSB -- Technical University of Ostrava, Ostrava, Czech Republic\\
\small Department of Algebra, Charles University, Prague, Czech Republic\\
\small \texttt{yifan.zhang@osu.cz}
}

\date{}

\begin{document}

\maketitle

\begin{abstract}
We study the minimisation of the positive principal eigenvalue for an indefinite-weight problem with asymmetric Robin parameters. The model is motivated by diffusive logistic equations in spatially heterogeneous environments, where the weight describes allocatable favourable resources and the Robin parameters measure boundary loss. After recalling the variational setting and the bang--bang reduction, we analyse the one-dimensional optimisation problem: the optimal favourable set is an interval, and the placement problem is reduced to a branchwise criterion. The key analytical tool is a shape-derivative formula for $a\mapsto\lambda(a)$, which shows that interior candidates are characterised by equality of the endpoint values of the positive eigenfunction, equivalently by the coupled transfer-matrix equations $f=0$ and $g=0$. We also introduce a Schrödinger-type extension with a fixed nonnegative background potential. In the coercive case we establish the corresponding principal-eigenvalue and bang--bang results, and in one dimension with constant potential we prove a compactness-type stability result showing that minimisers for small background potential converge, along subsequences, to minimisers of the unperturbed problem. No placement classification is claimed for general positive background potential. The computations are presented as numerical illustrations generated with an adaptive root-search protocol.
\end{abstract}

\medskip
\noindent\textbf{Keywords:} principal eigenvalue; indefinite weight; asymmetric Robin parameters; shape optimisation; transfer--matrix method; Schr\"odinger-type perturbation

\medskip
\noindent\textbf{2020 Mathematics Subject Classification:} 35P15, 35J25, 49Q10, 92D25.

\bigskip

\section{Introduction}\label{sec:introduction}

Spatial heterogeneity, diffusion, and boundary exchange are basic features of many reaction--diffusion models arising in population dynamics, spatial ecology, and related applied fields \cite{Skellam_1951,OkuboLevin2001,Murray2002,CantrellCosner2003,Hess1991}. In such models, dispersal is described by diffusion, local growth or decay depends on the position $x\in\Omega$, and the interaction with the exterior is encoded through boundary conditions. A standard prototype is the diffusive logistic equation on a bounded habitat,
\begin{align}\label{eq:logistic-intro}
\begin{cases}
\displaystyle u_t = \Delta u + \omega u[m(x)-u] & \text{in } \Omega\times \mathbb R_+,\\[0.4ex]
\partial_{\mathbf n}u + \beta_0 u = 0 & \text{on } \partial\Omega\times \mathbb R_+,\\[0.4ex]
u(x,0)\ge 0,\quad u(x,0)\not\equiv 0 & \text{in } \overline\Omega,
\end{cases}
\end{align}
where $\omega>0$ measures the overall reproductive intensity, $m$ is the intrinsic growth rate, and $\beta_0$ measures the inhospitableness of the boundary. When $\beta_0=0$ the boundary is impermeable, while positive values of $\beta_0$ model leakage, exchange, or mortality at the boundary; such Robin conditions are also standard in broader diffusion, thermal, and quantum models \cite{CantrellCosner2006,Daners2000,MunozCastaneda2015,Clarte2021}. In a heterogeneous habitat, the sign-changing coefficient $m$ distinguishes favourable and unfavourable regions, and the principal eigenvalue of the linearised operator determines the threshold between extinction and persistence \cite{Cantrell_1989,Smoller_1994,BerestyckiNirenbergVaradhan1994,Lou2006}.

This leads to the indefinite-weight eigenvalue problem
\begin{equation}\label{eq:pb}
\left\{
\begin{aligned}
\Delta \phi + \lambda m\phi &= 0 && \text{in } \Omega,\\
\partial_{\mathbf n}\phi + \beta\phi &= 0 && \text{on } \partial\Omega,
\end{aligned}
\right.
\end{equation}
where $\Omega\subset \mathbb R^n$ is a bounded Lipschitz domain, $\mathbf n$ is the unit outer normal on $\partial\Omega$, the weight $m\in L^\infty(\Omega)$ changes sign in $\Omega$, and
\[
-1\le m(x)\le \kappa \qquad\text{for a.e. }x\in \Omega,
\]
for some fixed constant $\kappa>0$. Throughout the paper, $\beta:\partial\Omega\to\mathbb R$ is a piecewise continuous Robin parameter. Following \cite{Brown1980,Senn1982,Afrouzi1999,Daners2013}, a real number $\lambda$ is called a \emph{principal eigenvalue} of \eqref{eq:pb} if the associated eigenfunction $\phi\in H^1(\Omega)$ is positive.

The optimisation problem appears when the total amount of favourable habitat is limited. One imposes pointwise bounds on $m$ together with an integral constraint on the total resource,
\[
\int_\Omega m \le -m_0 |\Omega|,
\]
with $m_0$ fixed, and then asks how the favourable part of the habitat should be placed so as to minimise the positive principal eigenvalue. From the modelling point of view, the amount of favourable resource is prescribed, but its spatial arrangement is free; lowering the principal eigenvalue lowers the persistence threshold and therefore enlarges the range of viable reproduction parameters. Questions of this type belong to the broader spectral-optimisation tradition going back at least to B\^ocher and Kre\u{\i}n \cite{Bocher1914,Krein1955}. In the indefinite-weight setting, the modern variational theory was developed for Neumann boundary conditions by Lou and Yanagida \cite{Lou2006}, for Robin boundary conditions by Hinterm\"uller, Kao and Laurain \cite{Hintermueller2011}, and further analysed by Lamboley, Laurain, Nadin and Privat \cite{Lamboley2016}. Beyond the interval case, complete structural results remain rare; the best understood higher-dimensional examples include cylindrical domains \cite{Kao2008} and spherical shells \cite{Schneider_2025}.

{The first aim of the present paper is to treat the one-dimensional problem with asymmetric Robin parameters, allowing distinct nonnegative parameters at the two endpoints. This setting is natural when the two ends of the habitat interact differently with the exterior. The main analytical issue is to determine the placement of the favourable interval. The proof follows the actual branch \(a\mapsto\lambda(a)\): a shape-derivative formula identifies the sign of \(\lambda' (a)\) with the difference of the squared eigenfunction values at the two ends of the favourable interval, and the transfer-matrix calculation converts the interior critical-point condition into the coupled equations \(f=0\) and \(g=0\).}
{The second aim is to record a natural Schr\"odinger-type extension in which the optimisable indefinite weight is accompanied by a fixed hostile background potential. In ecological terms, this additional potential may represent mortality, harvesting, pollution, predation, or some other adverse effect that is not itself subject to optimisation \cite{CuiLiMeiShi2017}. In the coercive case we show that the principal-eigenvalue theory and the bang--bang structure remain valid. We then study the one-dimensional family with constant potential and prove a compactness-type stability result showing that minimisers for small positive potential converge, along subsequences, to minimisers of the unperturbed problem. A second numerical study, again based on transfer matrices, illustrates this perturbative picture for representative canonical regimes.}

{We do not claim a placement classification for general positive background potentials; in the Schr\"odinger-type part, the placement statement is limited to the one-dimensional constant-potential compactness result.}

{To clarify the scope of the results, the general variational existence, activity, and bang--bang principles are recalled with the adaptations needed for the present notation and for the later Schr\"odinger-type extension. The new contribution of Section~\ref{sec:onedim} is the asymmetric-endpoint transfer-matrix analysis combined with a shape-derivative criterion along the actual eigenvalue branch. The computations in Section~\ref{sec:shooting} are numerical illustrations generated from the adaptive shooting code; they are not used as proofs of the analytical assertions.}

{The paper is organised as follows. After the variational preliminaries and the bang--bang reduction in the remainder of this section, Section~\ref{sec:onedim} contains the one-dimensional transfer-matrix analysis for the original indefinite-weight problem. Section~\ref{sec:potential-extension} introduces the Schr\"odinger-type extension, first at the general variational level and then for the one-dimensional constant-potential family. Section~\ref{sec:shooting} presents numerical illustrations: first for the unperturbed one-dimensional criteria, and then for the perturbative constant-potential regime. We conclude with a brief discussion of the analytical scope of the results, their modelling interpretation, and possible directions for further work.}

We first prove the existence and uniqueness of a positive principal eigenvalue following the method in \cite{Afrouzi1999}.

\begin{lemma}[\cite{Smoller_1994}, Theorem 11.4 and 11.10]\label{lem:smoller}
	For $\phi\in H^1(\Omega)$, let
	$$Q_\phi(\lambda):={\int_\Omega |\nabla \phi|^2  + \int_{\partial\Omega} \beta\phi^2 -\lambda \int_\Omega m\phi^2},$$
	and 
	$$\mu(\lambda):=\inf_{\phi\in H^1(\Omega)}\frac{ Q_\phi(\lambda)}{\int_\Omega \phi^2}.$$ Then
	$\lambda$ is a principal eigenvalue of \eqref{eq:pb}, if and only if $\mu(\lambda)=0$. 
\end{lemma}

\begin{lemma}[\cite{Afrouzi1999}, Lemma 2]\label{lem:mu-prime}
	If $\phi_0\in H^1(\Omega)$ is a minimiser of $Q_\phi(\lambda)$, i.e. $\mu(\lambda)=\frac{Q_{\phi_0}(\lambda)}{\int_\Omega \phi_0^2}$, then
	$$\mu'(\lambda)=-\frac{\int_\Omega m\phi_0^2}{\int_\Omega \phi_0^2}.$$
\end{lemma}

For any $\phi\in H^1(\Omega)$, the map $\lambda\mapsto Q_\phi(\lambda)$ is affine, hence concave. Since 
$$\mu(\lambda)=\inf\left\{Q_\phi(\lambda);\ \phi\in H^1(\Omega), \int_\Omega \phi^2=1\right\},$$
the map $\mu:\lambda\mapsto \mu(\lambda)$ is also concave as the infimum of a collection of concave functions. In other words, $\mu''(\lambda)\le 0$. Meanwhile, it is clear that there exist $\phi_1,\phi_2\in H^1(\Omega)$ satisfying $\int_\Omega m\phi_1^2>0$ and $\int_\Omega m\phi_2^2<0$, hence $\lim_{\lambda\to+\infty} \mu(\lambda)=\lim_{\lambda\to-\infty} \mu(\lambda)=-\infty$. Thus, $\mu=\mu(\lambda)$ increases until it attains its maximum, and decreases thereafter.

{When $\beta>0$ on a relatively measurable subset $\partial\Omega_\beta^+\subseteq \partial \Omega$ with positive surface measure, $\mathcal H^{n-1}(\partial\Omega_\beta^+)>0$, using the variational characterisation of the first Robin eigenvalue $\lambda_1(\beta)>0$, we have}
$$Q_\phi(0)\ge \lambda_1(\beta) \int_\Omega \phi^2>0$$
whenever $\int_\Omega \phi^2>0$,
hence $\mu(0)>0$, so $\mu$ has exactly two zeros, one positive and the other negative. In other words, $\eqref{eq:pb}$ has a unique positive principal eigenvalue and a unique negative principal eigenvalue. 

When $\beta\equiv 0$ a.e. on $\partial\Omega$, the situation is similar to Neumann boundary condition, and by taking the constant $\phi_c\equiv c$, we see 
$$Q_{\phi_c}(0)=0+c^2\int_{\partial\Omega}\beta = 0,$$
and thus $\mu(0)=0$. Since the associated eigenfunction is clearly a constant, by Lemma \ref{lem:mu-prime}, we get: 
\begin{itemize}
	\item if $\int_\Omega m<0$, then $\mu'(0)>0$, and $\mu$ has a unique positive zero corresponding to the unique positive principal eigenvalue of \eqref{eq:pb};
	\item if $\int_\Omega m>0$, then $\mu'(0)<0$, and $\mu$ has a unique negative zero corresponding to the unique negative principal eigenvalue of \eqref{eq:pb};
	\item if $\int_\Omega m=0$, then $\mu'(0)=0$, and $\mu$ has no other zeros, hence \eqref{eq:pb} has a unique principal eigenvalue $0$.
\end{itemize}

It is also well-known (cf. \cite{Brown1980}, \cite{Senn1982}, \cite{Lou2006} and \cite{Lamboley2016}) that the eigenvalue can be expressed with a Rayleigh quotient. To sum up, we have:

\begin{theorem}\label{thm:principal}
	If $\beta\ge 0$ on $\partial\Omega$ and {$\beta>0$ on a relatively measurable subset of $\partial\Omega$ of positive surface measure}, the problem \eqref{eq:pb} admits a unique positive principal eigenvalue $\lambda$, depending on $m$. Moreover, 
	$$\lambda(m) = \inf_{\phi\in \mathcal{S}(m)} \frac{\int_\Omega |\nabla \phi|^2 + \int_{\partial\Omega}\beta\phi^2}{\int_\Omega m\phi^2},$$
	where $\mathcal{S}(m):=\left\{\phi\in H^1(\Omega); \int_\Omega m\phi^2>0\right\}$.
\end{theorem}

Principal eigenvalues also admit monotonicity\footnote{We remark that an attempt to derive strict monotonicity in \cite{Afrouzi_2002}, Theorem 2.4, fails, for a positive eigenfunction could vanish on the boundary.} {with respect to the Robin parameter $\beta$.} 
\begin{theorem}\label{thm:monotone}
	Suppose $\lambda_1, \lambda_2$ are principal eigenvalues of problem \eqref{eq:pb} with Robin parameters $\beta_1$, $\beta_2$, respectively. {Denote by $\phi_1$ and $\phi_2$ the corresponding positive eigenfunctions.} Assuming $\int_\Omega m\phi_1^2>0$, if $\beta_1\ge \beta_2$ a.e. on $\partial\Omega$, then $\lambda_1\ge \lambda_2$. {If, in addition, there is a relatively measurable set $\Gamma\subseteq\partial\Omega$ of positive surface measure such that $\beta_1>\beta_2$ a.e. on $\Gamma$ and the trace of $\phi_1$ is positive on a subset of $\Gamma$ of positive surface measure, then $\lambda_1>\lambda_2$.} 
\end{theorem}
\begin{proof}
	By Lemma \ref{lem:smoller},
	\begin{align*}
		&\int_\Omega |\nabla\phi_1|^2+\int_{\partial\Omega}\beta_1\phi_1^2-\lambda_1\int_\Omega m\phi_1^2\\
		&=0\\
		&=\inf_{\phi\in H^1(\Omega)} \left(\int_\Omega |\nabla\phi|^2+\int_{\partial\Omega}\beta_2\phi^2-\lambda_2\int_\Omega m\phi^2\right)\\
		&\le \int_\Omega |\nabla\phi_1|^2+\int_{\partial\Omega}\beta_2\phi_1^2-\lambda_2\int_\Omega m\phi_1^2,
	\end{align*} 
	so 
	$$\lambda_1-\lambda_2\ge \frac{\int_{\partial\Omega} (\beta_1-\beta_2)\phi_1^2}{\int_\Omega m\phi_1^2}\ge 0.$$
	{Under the additional strictness assumption on $\Gamma$, the numerator is strictly positive, while the denominator is positive by hypothesis; hence $\lambda_1> \lambda_2$.}
\end{proof}

We are interested in the minimisation of the positive principal eigenvalue $\lambda(m)$. {Here and below we use the notation $\Omega_m^+:=\{x\in\Omega:\ m(x)>0\}$.}
\begin{equation}\label{eq:pb-v1}
	\begin{aligned}
		\inf_{m\in \mathcal{M}_{m_0,\kappa}} &\lambda(m),\\
		\mathcal{M}_{m_0,\kappa} &:= \left\{m\in L^\infty(\Omega);\ -1\le m\le \kappa, \left|\Omega_m^+\right|>0, \int_\Omega m\le -m_0|\Omega|\right\}
	\end{aligned}
\end{equation}
Here $m_0$ is a real constant such that $m_0\in(-\kappa, 1)$ {if $\beta$ is positive on a subset of $\partial\Omega$ of positive surface measure} and $m_0\in(0,1)$ {in the pure Neumann case $\beta\equiv0$}. 

{We recall the Robin version of the standard existence and bang--bang result, and include the argument because we use the quantifier below.} {The stronger ``every minimiser'' form is the one proved in the Neumann case in \cite[Theorem~1.1]{Lou2006} and stated for the Robin case in \cite[Theorem~1]{Hintermueller2011}.}
\begin{theorem}\label{thm:bangbang}
	{The infimum in \eqref{eq:pb-v1} is attained. Moreover, every minimiser is bang--bang: if $m^*\in \mathcal{M}_{m_0,\kappa}$ satisfies $\lambda(m^*)=\inf_{m\in\mathcal{M}_{m_0,\kappa}}\lambda(m)$, then there exists a measurable subset $E^*$ of $\Omega$ such that}
	\begin{equation}
		m^* = \kappa \chi_{E^*} - \chi_{\Omega\setminus E^*} \text{ a.e. } x\in \Omega,
	\end{equation}
	where $\chi_U$ denotes the characteristic function of a set $U$. {Every minimiser also satisfies the active volume constraint}
	$$\int_\Omega m^* = -m_0|\Omega|.$$
\end{theorem}

\begin{proof}
{Existence is the standard direct-method argument for indefinite weights; the fixed nonnegative Robin term does not change the compactness and lower-semicontinuity argument used in \cite{Lou2006,Hintermueller2011}. We prove the stronger bang--bang assertion for an arbitrary minimiser. Let $\widehat m\in\mathcal M_{m_0,\kappa}$ be any minimiser and let $u>0$ be an associated eigenfunction. Set
\[
A(u):=\int_\Omega |\nabla u|^2+\int_{\partial\Omega}\beta u^2,
\qquad
D_m(u):=\int_\Omega m u^2.
\]
Then $D_{\widehat m}(u)>0$ and $\lambda(\widehat m)=A(u)/D_{\widehat m}(u)$. For every admissible competitor $m$ with $D_m(u)>0$, the fixed-eigenfunction comparison gives
\[
\lambda(\widehat m)
\le \lambda(m)
\le \frac{A(u)}{D_m(u)}.
\]
Hence $D_m(u)\le D_{\widehat m}(u)$. Thus $\widehat m$ maximises the linear functional $m\mapsto D_m(u)$ among admissible competitors with positive denominator.

We first prove activity. If $\int_\Omega \widehat m<-m_0|\Omega|$, then, since $\widehat m\le\kappa$ and $m_0>-\kappa$, there is a measurable subset of positive measure on which $\widehat m<\kappa$. Increasing $\widehat m$ by a sufficiently small positive amount on such a subset gives an admissible competitor still having positive denominator, but with strictly larger value of $D_m(u)$ because $u>0$ a.e. This contradicts the maximality just proved. Therefore
\[
\int_\Omega \widehat m=-m_0|\Omega|.
\]

Now apply the bathtub principle to the linear functional $m\mapsto D_m(u)$ under the bounds $-1\le m\le\kappa$ and the active mass constraint. There exists a maximiser $m^\sharp$ of the form
\[
m^\sharp=\kappa\chi_{E^\sharp}-\chi_{\Omega\setminus E^\sharp},
\qquad |E^\sharp|=c|\Omega|.
\]
Moreover $D_{m^\sharp}(u)\ge D_{\widehat m}(u)>0$. Therefore
\[
\lambda(\widehat m)
\le \lambda(m^\sharp)
\le \frac{A(u)}{D_{m^\sharp}(u)}
\le \frac{A(u)}{D_{\widehat m}(u)}
=\lambda(\widehat m).
\]
All inequalities are equalities. In particular, $m^\sharp$ is also a minimiser and $u$ attains the Rayleigh quotient for $m^\sharp$, so $u$ is an eigenfunction for $m^\sharp$ with the same eigenvalue. Hence, in the weak sense,
\[
\Delta u+\lambda(\widehat m)\widehat m u=0,
\qquad
\Delta u+\lambda(\widehat m)m^\sharp u=0.
\]
Subtracting gives $(m^\sharp-\widehat m)u=0$ a.e. in $\Omega$. Since the principal eigenfunction is positive a.e. in $\Omega$, we obtain $\widehat m=m^\sharp$ a.e. Thus the arbitrary minimiser $\widehat m$ is bang--bang and satisfies the active constraint.}
\end{proof}

{In other words, every optimal weight function has a bang--bang distribution in $\Omega$ with respect to some measurable optimal set $E^*$.} The problem can therefore be formulated as a shape optimisation one:
\begin{equation}\label{eq:pb-v2}
	\begin{aligned}
		\inf_{E\in \mathcal{E}_{c,\kappa}}  \lambda(E) & :=\lambda(\kappa \chi_E - \chi_{\Omega\setminus E}),\\
		\mathcal{E}_{c,\kappa} &:= \{E\subseteq \Omega; E\ \text{is measurable}, |E|=c|\Omega|\}
	\end{aligned}
\end{equation}
where $c:=\frac{1-m_0}{1+\kappa}$ and $c\in (0,1)$ {if $\beta$ is positive on a subset of $\partial\Omega$ of positive surface measure}, while $c\in\left(0,\frac{1}{1+\kappa}\right)$ {in the pure Neumann case $\beta\equiv0$}. {The equality $|E|=c|\Omega|$ is justified by the activity of the resource constraint in Theorem~\ref{thm:bangbang}.} The optimal $E$ is expressed as a function of $\beta$, $\kappa$ and $c$. 

In \cite{Hintermueller2011} and \cite{Lamboley2016}, this minimisation problem, formulated as a shape optimisation one, is thoroughly investigated {under a constant Robin parameter}, with the one-dimensional case fully solved. Unfortunately, for higher dimensional cases, few results are known. In \cite{Hintermueller2011}, they studied this problem in a cylindrical domain with Robin boundary condition on the top and bottom and Neumann boundary condition on the lateral boundary. In \cite{Lamboley2016}, it is shown that when $\partial\Omega$ is connected, the optimal set is not always a ball. In \cite{Schneider_2025}, the 1-dimensional result is generalised onto spherical shells of any dimension. 

\section{The 1-dimensional case}\label{sec:onedim}

{In this section, we analyse problem \eqref{eq:pb} in one dimension. After the interval reduction, the admissible favourable set is written as \(E_a=(a,a+c)\), and the remaining problem is to minimise the principal eigenvalue along the branch \(a\mapsto\lambda(a)\). The main point of the argument is that all sign information is taken along this branch, not at a frozen value of \(\lambda\).} 

Consider the problem \eqref{eq:pb} in $\Omega=(0,1)\subseteq \mathbb R$. {The following one-dimensional interval reduction is the only global structural input used in this section.} {In the proof we write \(\beta_0,\beta_1\) for the two endpoint Robin parameters, so that the boundary conditions are \(u'(0)-\beta_0u(0)=0\) and \(u'(1)+\beta_1u(1)=0\), with \(\beta_0,\beta_1\ge0\).}

\begin{proposition}[Interval reduction]\label{prop:optimal-interval}
	{If $\Omega=(0,1)$ and the endpoint Robin parameters satisfy $\beta_0,\beta_1\ge0$, then every optimal set $E^*$ for \eqref{eq:pb-v2} is, up to a set of Lebesgue measure zero, an interval of length $c$.}
\end{proposition}

\begin{proof}
{We adapt the two-sided monotone rearrangement argument used in the proof of \cite[Proposition~4]{Lamboley2016}. We give the details because the Robin coefficients at the two endpoints are now allowed to be different. By Theorem~\ref{thm:bangbang}, every minimising weight is bang--bang. Hence, for an arbitrary optimal set $E\subset(0,1)$ with $|E|=c$, we may write}
\[
{m_E=\kappa\chi_E-\chi_{(0,1)\setminus E}.}
\]
{Let $u>0$ be the associated principal eigenfunction. We use the Rayleigh quotient}
\[
{\mathcal R_{\beta_0,\beta_1}(v,E):=
\frac{\displaystyle\int_0^1 |v'|^2+\beta_0v(0)^2+\beta_1v(1)^2}
{\displaystyle\int_0^1 m_Ev^2},}
\]
{for test functions for which the denominator is positive. Since $u$ is the positive principal eigenfunction for $E$,}
\[
{D_E(u):=\int_0^1m_Eu^2=(\kappa+1)\int_Eu^2-\int_0^1u^2>0,}
\]
{and $\lambda(E)=\mathcal R_{\beta_0,\beta_1}(u,E)$.}

{We first prove that this same set $E$ is selected by the bathtub principle for the fixed density $u^2$. Keeping $u$ fixed, the numerator of $\mathcal R_{\beta_0,\beta_1}(u,F)$ is independent of the admissible set $F$, while}
\[
{D_F(u)=(\kappa+1)\int_Fu^2-\int_0^1u^2.}
\]
{If there were a measurable set $F\subset(0,1)$ with $|F|=c$ and $\int_Fu^2>\int_Eu^2$, then $D_F(u)>D_E(u)>0$ and therefore}
\[
{\lambda(F)\le \mathcal R_{\beta_0,\beta_1}(u,F)<\mathcal R_{\beta_0,\beta_1}(u,E)=\lambda(E),}
\]
{contradicting the optimality of $E$. Hence $E$ maximises $F\mapsto\int_Fu^2$ among all measurable sets of measure $c$. By the bathtub principle, there are $t\ge0$ and a measurable set $G\subseteq\{u=t\}$ such that}
\[
{E=\{u>t\}\cup G,
\qquad |E|=c.}
\]
{This representation applies to the original arbitrary optimiser $E$, up to null sets. Moreover, the level set $\{u=t\}$ has zero measure. Indeed, $u$ satisfies a one-dimensional equation with bounded coefficient, hence $u\in W^{2,\infty}(0,1)$. If $\{u=t\}$ had positive measure, then $u'=0$ and $u''=0$ a.e. on that level set. The equation $u''+\lambda m_Eu=0$, together with $\lambda>0$, $u>0$ and $m_E\in\{\kappa,-1\}$, gives a contradiction. Thus $E=\{u>t\}$ up to a null set.}

{Let $x_0\in[0,1]$ be a point where $u$ attains its maximum. On $(0,x_0)$ we take the increasing rearrangement of $u$, and on $(x_0,1)$ we take the decreasing rearrangement of $u$; if $x_0=0$ or $x_0=1$, the empty side is omitted. Denote the resulting two-sided monotone rearrangement by $u^\sharp$. Then $u^\sharp$ is equimeasurable with $u$, belongs to $H^1(0,1)$, is nondecreasing on $(0,x_0)$ and nonincreasing on $(x_0,1)$, and the one-dimensional P\'olya inequality applied on the two sides gives}
\[
{\int_0^1 |(u^\sharp)'|^2\le \int_0^1 |u'|^2.}
\]
{The endpoint traces are also controlled separately: the increasing rearrangement on the left places the smallest values at $0$, and the decreasing rearrangement on the right places the smallest values at $1$. In the endpoint-maximum cases the corresponding side is absent and the trace is unchanged. Hence}
\[
{u^\sharp(0)^2\le u(0)^2,
\qquad
u^\sharp(1)^2\le u(1)^2.}
\]
{Since $\beta_0,\beta_1\ge0$,}
\[
{\beta_0u^\sharp(0)^2+\beta_1u^\sharp(1)^2
\le
\beta_0u(0)^2+\beta_1u(1)^2.}
\]
{This is the only point at which the Robin boundary term enters the rearrangement comparison. The estimate is componentwise at the two endpoints, so equality of the Robin coefficients is not required.}

{Because $u^\sharp$ is equimeasurable with $u$, the maximum value of $F\mapsto\int_F(u^\sharp)^2$ over all sets $F$ of measure $c$ is the same as the corresponding maximum for $u$. Since $u^\sharp$ is unimodal, this maximum is attained by a superlevel set of $u^\sharp$, which is an interval; denote it by $E^\sharp$. Then $|E^\sharp|=c$ and}
\[
{\int_{E^\sharp}(u^\sharp)^2=\int_Eu^2,
\qquad
\int_0^1 (u^\sharp)^2=\int_0^1u^2,}
\]
{so that $D_{E^\sharp}(u^\sharp)=D_E(u)>0$. Combining this equality of denominators with the Dirichlet-energy and boundary-term comparisons yields}
\[
{\lambda(E^\sharp)
\le \mathcal R_{\beta_0,\beta_1}(u^\sharp,E^\sharp)
\le \mathcal R_{\beta_0,\beta_1}(u,E)
=\lambda(E).}
\]
{Since $E$ is optimal, equality holds throughout. In particular, equality holds in the one-dimensional P\'olya inequality on both sides of $x_0$. By the equality case used in \cite[Proposition~4]{Lamboley2016}, the original eigenfunction $u$ is itself nondecreasing on $(0,x_0)$ and nonincreasing on $(x_0,1)$, after modification on a null set; by continuity of $u$, this monotonicity holds for the continuous representative. Therefore the superlevel set $\{u>t\}$ is an interval, possibly with one endpoint degenerated. Since $E=\{u>t\}$ up to a null set and $|E|=c$, the arbitrary optimal set $E$ is an interval of length $c$ up to a null set.}
\end{proof}

{The formulae below are first derived for $0<a<a+c<1$. Endpoint placements $a=0$ and $a=1-c$ are included by the same transfer matrices with one zero-length exterior interval, or equivalently by continuity of $a\mapsto\lambda(a)$ established below.}
Now let the optimal set be the interval $E=(a,b)$ where $0<a<b<1$ and $c=b-a$. We determine the optimal value of $a$ by adapting the method of \cite{Hintermueller2011}. Let $\beta_0=\beta(0)\ge 0$ and {$\beta_1=\beta(1)\ge0$, with $\beta_0+\beta_1>0$}. Rewrite the problem using $m=\kappa\chi_{E}-\chi_{\Omega\setminus E}$, and we have

\begin{equation}\label{eq:pb-1d}
	\left\{
	\begin{aligned}
		\phi''-\lambda\phi &= 0&\ \text{in} &\ (0,a), \\
		\phi''+\lambda\kappa\phi&=0&\ \text{in} &\ (a,b),\\
		\phi''-\lambda\phi &= 0&\ \text{in} &\ (b,1), \\
		\phi'(0)-\beta_0\phi(0)&=0,\\
		\phi'(1)+\beta_1\phi(1)&=0,\\
		\phi(a)^+-\phi(a)^-&=0,\\
		\phi(b)^+-\phi(b)^-&=0,\\
		\phi'(a)^+-\phi'(a)^-&=0,\\
		\phi'(b)^+-\phi'(b)^-&=0.
	\end{aligned}
	\right.
\end{equation}

Since $\lambda>0$, $\phi>0$, we may assume
\begin{equation}\label{eq:sol-1d}
	\phi(x)=	\left\{
	\begin{aligned}
		C_1 \cosh\left(\sqrt{\lambda}(x-a)\right) + C_2 \sinh\left(\sqrt{\lambda}(x-a)\right) \ \text{in}\ (0,a),\\
		C_3 \cos\left(\sqrt{\lambda\kappa}(x-a)\right) + C_4\sin\left(\sqrt{\lambda\kappa}(x-a)\right)\ \text{in}\ (a,b),\\
		C_5\cosh\left(\sqrt{\lambda}(x-b)\right)+C_6\sinh\left(\sqrt{\lambda}(x-b)\right)\ \text{in}\ (b,1),
	\end{aligned}
	\right.
\end{equation}
where $C_i=C_i(E)$ are constants for $1\le i\le 6$. We now simplify the boundary and transmission conditions in \eqref{eq:pb-1d}. Collecting the coefficients of $C_i$ ($1\le i\le 6$), we obtain
\small
\begin{equation}
	\left\{
	\begin{aligned}
		&\left(-\sqrt{\lambda}\sinh(a\sqrt{\lambda})-\beta_0\cosh(a\sqrt{\lambda})\right)C_1 + \left(\sqrt{\lambda}\cosh(a\sqrt{\lambda})+\beta_0\sinh(a\sqrt{\lambda})\right)C_2=0\\
		&\left(\sqrt{\lambda}\sinh(\sqrt{\lambda}(1-b))+\beta_1\cosh(\sqrt{\lambda}(1-b))\right)C_5\\
		&\ \ \ +\left(\sqrt{\lambda}\cosh(\sqrt{\lambda}(1-b))+\beta_1\sinh(\sqrt{\lambda}(1-b))\right)C_6 =0\\
		&C_3-C_1=0\\
		&C_5-C_3\cos(\sqrt{\lambda\kappa}c) - C_4\sin(\sqrt{\lambda\kappa}c)=0\\
		&\sqrt{\kappa}C_4-C_2=0\\
		&C_6+C_3\sqrt{\kappa}\sin(\sqrt{\lambda\kappa}c)-C_4\sqrt{\kappa}\cos(\sqrt{\lambda\kappa}c)=0
	\end{aligned}
	\right.
\end{equation}
\normalsize
Using the last 4 equations, we can express $C_1,C_2,C_5,C_6$ in terms of $C_3$ and $C_4$, and then we obtain two linear equations of $C_3$ and $C_4$ only:
\small
\begin{equation}
	\left\{
	\begin{aligned}
		&\left(-\sqrt{\lambda}{\mathfrak{s}}_\lambda(a)-\beta_0{\mathfrak{c}}_\lambda(a)\right)C_3 + \sqrt{\kappa}\left(\sqrt{\lambda}{\mathfrak{c}}_\lambda(a)+\beta_0{\mathfrak{s}}_\lambda(a)\right)C_4=0\\
		&\left[\cos(\sqrt{\lambda\kappa}c)\left(\sqrt{\lambda}{\mathfrak{s}}_\lambda(b)+\beta_1{\mathfrak{c}}_\lambda(b)\right)-\sqrt{\kappa}\sin(\sqrt{\lambda\kappa}c)\left(\sqrt{\lambda}{\mathfrak{c}}_\lambda(b)+\beta_1{\mathfrak{s}}_\lambda(b)\right)\right]C_3\\
		&+\left[\sin(\sqrt{\lambda\kappa}c)\left(\sqrt{\lambda}{\mathfrak{s}}_\lambda(b)+\beta_1{\mathfrak{c}}_\lambda(b)\right)+\sqrt{\kappa}\cos(\sqrt{\lambda\kappa}c)\left(\sqrt{\lambda}{\mathfrak{c}}_\lambda(b)+\beta_1{\mathfrak{s}}_\lambda(b)\right)\right]C_4=0
	\end{aligned}
	\right.
\end{equation}
\normalsize
where ${\mathfrak{s}}_\lambda(a):=\sinh(\sqrt{\lambda}a)$, ${\mathfrak{c}}_\lambda(a):=\cosh(\sqrt{\lambda}a)$,  ${\mathfrak s}_\lambda(b):=\sinh(\sqrt{\lambda}(1-b))$ and ${\mathfrak c}_\lambda(b):=\cosh(\sqrt{\lambda}(1-b))$. 
Considered as a $2\times 2$ matrix applied to the column vector $(C_3,C_4)^t$, they have a non-trivial solution if and only if the determinant of the matrix vanishes, i.e. 
\small
\begin{multline}\label{eq:det0-original}
	\left(-\sqrt{\lambda}{\mathfrak{s}}_\lambda(a)-\beta_0{\mathfrak{c}}_\lambda(a)\right)\left[\sin(\sqrt{\lambda\kappa}c)\left(\sqrt{\lambda}{\mathfrak{s}}_\lambda(b)+\beta_1{\mathfrak{c}}_\lambda(b)\right)\right.\\
	+\left.\sqrt{\kappa}\cos(\sqrt{\lambda\kappa}c)\left(\sqrt{\lambda}{\mathfrak{c}}_\lambda(b)+\beta_1{\mathfrak{s}}_\lambda(b)\right)\right]\\
	= \sqrt{\kappa}\left(\sqrt{\lambda}{\mathfrak{c}}_\lambda(a)+\beta_0\sinh(a\sqrt{\lambda})\right) \left[\cos(\sqrt{\lambda\kappa}c)\left(\sqrt{\lambda}{\mathfrak{s}}_\lambda(b)+\beta_1{\mathfrak{c}}_\lambda(b)\right)\right.\\
	-\left.\sqrt{\kappa}\sin(\sqrt{\lambda\kappa}c)\left(\sqrt{\lambda}{\mathfrak{c}}_\lambda(b)+\beta_1{\mathfrak{s}}_\lambda(b)\right)\right].
\end{multline}
\normalsize
Rearranging the terms gives
\small
\begin{multline}\label{eq:det0-ordered}
	\cosh(\sqrt{\lambda}(1-b))\left[\cosh(\sqrt{\lambda}a)\left(-(\beta_0+\beta_1)\sqrt{\lambda\kappa}\cos(\sqrt{\lambda\kappa}c)+(\lambda \kappa -\beta_0\beta_1)\sin(\sqrt{\lambda\kappa}c)\right) \right.\\
	\left. +\sinh(\sqrt{\lambda}a)\left(-\sqrt{\kappa}(\lambda+\beta_0\beta_1)\cos(\sqrt{\lambda\kappa}c)+(\kappa\beta_0-\beta_1)\sqrt{\lambda}\sin(\sqrt{\lambda\kappa}c)\right)\right]\\
	+\sinh(\sqrt{\lambda}(1-b))\left[\cosh(\sqrt{\lambda}a)\left(-\sqrt{\kappa}(\lambda+\beta_0\beta_1)\cos(\sqrt{\lambda\kappa}c)+(\kappa\beta_1-\beta_0)\sqrt{\lambda}\sin(\sqrt{\lambda\kappa}c)\right)\right.\\
	\left.+ \sinh(\sqrt{\lambda}a)\left(-\sqrt{\lambda\kappa}(\beta_0+\beta_1)\cos(\sqrt{\lambda\kappa}c)+(\beta_0\beta_1\kappa-\lambda)\sin(\sqrt{\lambda\kappa}c)\right)	\right]=0.
\end{multline}
\normalsize
It is useful to check \eqref{eq:det0-ordered} in two extreme cases. When $\beta_0=\beta_1=0$, the problem reduces to the Neumann case, and we obtain the characteristic equation\footnote{We remark that the characteristic equations (18) and (20) obtained in \cite{Hintermueller2011} based on {a constant Robin parameter} admit the same typo on the left-hand side. }
\begin{align}\label{eq:char-neumann}
	\tanh(\lambda^{\frac{1}{2}}(1-b))=\frac{\kappa^{\frac{1}{2}}\tan(\lambda^{\frac{1}{2}}\kappa^{\frac{1}{2}}c)-\tanh(\lambda^{\frac{1}{2}}a)}{1+\kappa^{-\frac{1}{2}}\tanh(\lambda^{\frac{1}{2}}a)\tan(\lambda^{\frac{1}{2}}\kappa^{\frac{1}{2}}c)},
\end{align}
which reduces, when $a=0$, to the equation (3.3) in \cite{Lou2006}:
\begin{align}
	\kappa^{\frac{1}{2}}\tan(\lambda^{\frac{1}{2}}\kappa^{\frac{1}{2}}c)=\tanh(\lambda^{\frac{1}{2}}(1-c)).
\end{align}
On the other hand, when $\beta_0,\beta_1\to+\infty$, the problem reduces to the Dirichlet case, and we similarly obtain the characteristic equation
\begin{align}\label{eq:char-dirichlet}
	\tanh(\lambda^{\frac{1}{2}}(1-b))=\frac{\kappa^{-\frac{1}{2}}\tan(\lambda^{\frac{1}{2}}\kappa^{\frac{1}{2}}c)+\tanh(\lambda^{\frac{1}{2}}a)}{\kappa^{\frac{1}{2}}\tanh(\lambda^{\frac{1}{2}}a)\tan(\lambda^{\frac{1}{2}}\kappa^{\frac{1}{2}}c)-1},
\end{align}
which reduces, when $a=0$, to the equation (21) in \cite{Hintermueller2011}:
\begin{align}
	\tan(\lambda^{\frac{1}{2}}\kappa^{\frac{1}{2}}c)=-\kappa^{\frac{1}{2}}\tanh(\lambda^{\frac{1}{2}}(1-c)).
\end{align}

We next use the product-to-sum formulae for hyperbolic functions to simplify this expression. Using $a\sqrt{\lambda} + \sqrt{\lambda}(1-b) = \sqrt{\lambda}(1-c)$ we get
\small
\begin{multline}\label{eq:det0-pdtosum}
	0=(\kappa+1)(\lambda-\beta_0\beta_1)\cosh(\sqrt{\lambda}(a+b-1))\sin(\sqrt{\lambda\kappa}c)\\
	+ (\kappa+1)(\beta_0-\beta_1)\sqrt{\lambda}\sinh(\sqrt{\lambda}(a+b-1))\sin(\sqrt{\lambda\kappa}c)\\
	+\cosh(\sqrt{\lambda}(1-c))\left((\kappa-1)(\lambda+\beta_0\beta_1)\sin(\sqrt{\lambda\kappa}c)-2\sqrt{\lambda\kappa}(\beta_0+\beta_1)\cos(\sqrt{\lambda\kappa}c)\right)\\
	+\sinh(\sqrt{\lambda}(1-c))\left((\kappa-1)(\beta_0+\beta_1)\sqrt{\lambda}\sin(\sqrt{\lambda\kappa}c)-2\sqrt{\kappa}(\beta_0\beta_1+\lambda)\cos(\sqrt{\lambda\kappa}c)\right).
\end{multline}
\normalsize
We express the right-hand side of \eqref{eq:det0-pdtosum} as a function $f(a, \beta_0,\beta_1,\kappa,\lambda)$, using $b=a+c$:
\small
\begin{multline}\label{eq:f-vanish}
	f(a, \beta_0,\beta_1,\kappa,\lambda)=(\kappa+1)(\lambda-\beta_0\beta_1)\cosh(\sqrt{\lambda}(2a+c-1))\sin(\sqrt{\lambda\kappa}c)\\
	+ (\kappa+1)(\beta_0-\beta_1)\sqrt{\lambda}\sinh(\sqrt{\lambda}(2a+c-1))\sin(\sqrt{\lambda\kappa}c)\\
	+\cosh(\sqrt{\lambda}(1-c))\left((\kappa-1)(\lambda+\beta_0\beta_1)\sin(\sqrt{\lambda\kappa}c)-2\sqrt{\lambda\kappa}(\beta_0+\beta_1)\cos(\sqrt{\lambda\kappa}c)\right)\\
	+\sinh(\sqrt{\lambda}(1-c))\left((\kappa-1)(\beta_0+\beta_1)\sqrt{\lambda}\sin(\sqrt{\lambda\kappa}c)-2\sqrt{\kappa}(\beta_0\beta_1+\lambda)\cos(\sqrt{\lambda\kappa}c)\right).
\end{multline}
\normalsize
Taking partial derivative with respect to $a$ gives
\begin{multline}\label{eq:df-a}
	\partial_a f(a, \beta_0,\beta_1,\kappa,\lambda)=2(\kappa+1)\sqrt{\lambda}\sin(\sqrt{\lambda\kappa}c)\cosh(\sqrt{\lambda}(2a+c-1))\\
	\cdot\left[(\lambda-\beta_0\beta_1)\tanh\left(2\sqrt{\lambda}\left(a-\frac{1-c}{2}\right)\right)+\sqrt{\lambda}(\beta_0-\beta_1)\right].
\end{multline}

\begingroup\color{black}
We shall not use an a priori fixed spectral window such as
\(0<\lambda<\pi^2/(4c^2\kappa)\). Instead, the sign analysis below is carried out along the
actual branch \(a\mapsto\lambda(a)\). The formula \eqref{eq:df-a} remains useful,
but only after it is coupled with the branch \(\lambda=\lambda(a)\).

For later reference, set
\begin{align}\label{eq:df-a-g}
    g(a,\beta_0,\beta_1,\lambda):=(\lambda-\beta_0\beta_1)
    \tanh\left(2\sqrt{\lambda}\left(a-\frac{1-c}{2}\right)\right)
    +\sqrt{\lambda}(\beta_0-\beta_1).
\end{align}
Thus \eqref{eq:df-a} can be written as
\begin{align}\label{eq:df-a-factorised}
    \partial_a f(a,\beta_0,\beta_1,\kappa,\lambda)
    =2(\kappa+1)\sqrt{\lambda}\sin(\sqrt{\lambda\kappa}c)
    \cosh(\sqrt{\lambda}(2a+c-1))\,g(a,\beta_0,\beta_1,\lambda).
\end{align}

We now make the dependence on the placement explicit. For \(a\in[0,1-c]\), let
\begin{align}\label{eq:ma-def}
    m_a:=-1+(\kappa+1)\chi_{(a,a+c)}
    =\kappa\chi_{(a,a+c)}-\chi_{(0,1)\setminus(a,a+c)},
\end{align}
and denote by \(\lambda(a)\) the positive principal eigenvalue of \eqref{eq:pb-1d}
with \(m=m_a\). Let \(u_a>0\) be a corresponding eigenfunction.

\begin{lemma}\label{lem:spectral-positive-window}
For every \(a\in[0,1-c]\), the positive principal eigenvalue satisfies
\begin{align}\label{eq:true-positive-window}
    0<\sqrt{\lambda(a)\kappa}\,c<\pi .
\end{align}
In particular \(\sin(\sqrt{\lambda(a)\kappa}c)>0\).
\end{lemma}

\begin{proof}
{Let \(u_a>0\) be the principal eigenfunction and set \(\omega=\sqrt{\lambda(a)\kappa}\).
On the favourable interval \((a,a+c)\), the equation is \(u_a''+\omega^2u_a=0\). A
non-trivial solution of this equation can have a fixed strict sign only on intervals of
length strictly smaller than \(\pi/\omega\): after writing it as
\(A\cos(\omega x)+B\sin(\omega x)\), its consecutive zeros are separated by exactly
\(\pi/\omega\), and the equality case would force zeros at both endpoints of the sign
interval. Since the principal eigenfunction is strictly positive on the closed interval
\([a,a+c]\), we must have \(c<\pi/\omega\). Hence \(\omega c<\pi\). Since
\(\lambda(a)>0\) and \(\kappa>0\), the left-hand inequality is immediate.}
\end{proof}

\begin{lemma}\label{lem:lambda-regularity}
The map \(a\mapsto\lambda(a)\) is continuous on \([0,1-c]\) and real analytic on
\((0,1-c)\). Moreover, for \(a\in(0,1-c)\), the zero \(\lambda=\lambda(a)\) of the
characteristic equation \(f(a,\beta_0,\beta_1,\kappa,\lambda)=0\) is simple.
\end{lemma}

\begin{proof}
It is convenient to use the transfer residual rather than only the expanded expression
\eqref{eq:f-vanish}. Let
\[
P_-(L,\lambda):=
\begin{pmatrix}
\cosh(\sqrt\lambda L)&\sinh(\sqrt\lambda L)/\sqrt\lambda\\
\sqrt\lambda\sinh(\sqrt\lambda L)&\cosh(\sqrt\lambda L)
\end{pmatrix}
\]
for an unfavourable interval of length \(L\), and
\[
P_+(L,\lambda):=
\begin{pmatrix}
\cos(\sqrt{\kappa\lambda} L)&\sin(\sqrt{\kappa\lambda} L)/\sqrt{\kappa\lambda}\\
-\sqrt{\kappa\lambda}\sin(\sqrt{\kappa\lambda} L)&\cos(\sqrt{\kappa\lambda} L)
\end{pmatrix}
\]
for a favourable interval of length \(L\). These matrices are real analytic in
\((L,\lambda)\) for \(L\ge0\), \(\lambda>0\), with the usual limiting interpretation at
\(L=0\). The boundary residual
\[
\begin{aligned}
R(a,\lambda):={}&e_2P_-(1-a-c,\lambda)P_+(c,\lambda)P_-(a,\lambda)(1,\beta_0)^\top\\
&+\beta_1e_1P_-(1-a-c,\lambda)P_+(c,\lambda)P_-(a,\lambda)(1,\beta_0)^\top
\end{aligned}
\]
is therefore real analytic in \((a,\lambda)\) on \((0,1-c)\times(0,\infty)\), and its
zeros are precisely the eigenvalues of the regular one-dimensional separated
boundary-value problem.

For each fixed \(a\), the principal eigenvalue is a simple zero of this residual. Indeed,
let \(y(x,\lambda)\) be the solution with \(y(0,\lambda)=1\) and
\(y'(0,\lambda)=\beta_0\). At an eigenvalue \(\lambda=\lambda(a)\), the residual is
\(R(a,\lambda)=y'(1,\lambda)+\beta_1y(1,\lambda)\). Differentiating the equation
\(y''+\lambda m_a y=0\) with respect to \(\lambda\), and writing
\(z=\partial_\lambda y\), gives \(z''+\lambda m_a z=-m_a y\), with
\(z(0)=z'(0)=0\). The Lagrange identity yields
\[
\frac{d}{dx}\bigl(z'y-zy'\bigr)=-m_a y^2.
\]
After integration over \((0,1)\), and using the right Robin condition
\(y'(1)+\beta_1y(1)=0\), we obtain
\[
 y(1)\,\partial_\lambda R(a,\lambda(a))=-\int_0^1 m_a y^2\,dx.
\]
For the positive principal eigenvalue the integral on the right is positive, while
\(y(1)>0\). Hence \(\partial_\lambda R(a,\lambda(a))<0\). The principal zero is
therefore simple, and the implicit-function theorem yields real analyticity of
\(a\mapsto\lambda(a)\) on \((0,1-c)\).

Continuity at \(a=0\) and \(a=1-c\) follows from the same residual formula with a
zero-length exterior interval, or equivalently from the variational characterisation and
\(m_a\to m_{a_0}\) in \(L^1(0,1)\). {A direct comparison of the transfer residual above with the expanded determinant in \eqref{eq:f-vanish} gives
\[
    f(a,\beta_0,\beta_1,\kappa,\lambda)=-2\sqrt{\kappa\lambda}\,R(a,\lambda).
\]
The factor \(-2\sqrt{\kappa\lambda}\) is nonzero for \(\kappa>0\) and \(\lambda>0\). Therefore \(f\) and \(R\) have exactly the same positive zeros. At the principal branch,
\[
    \partial_\lambda f(a,\beta_0,\beta_1,\kappa,\lambda(a))
    =-2\sqrt{\kappa\lambda(a)}\partial_\lambda R(a,\lambda(a))\ne0,
\]
because the term obtained by differentiating the prefactor is multiplied by \(R(a,\lambda(a))=0\). Thus the expanded equation \(f=0\) has the same simple principal zero as the residual equation.}
\end{proof}

\begin{proposition}[Shape derivative]\label{prop:shape-derivative}
Fix \(a\in(0,1-c)\), and normalise the positive eigenfunction by
\begin{align}\label{eq:normalization-D}
    \int_0^1 m_a u_a^2\,dx=1.
\end{align}
Then
\begin{align}\label{eq:shape-derivative}
    \lambda'(a)=\lambda(a)(\kappa+1)\big(u_a(a)^2-u_a(a+c)^2\big).
\end{align}
Consequently, an interior critical point satisfies
\begin{align}\label{eq:critical-endpoint-values}
    \lambda'(a)=0 \quad\Longleftrightarrow\quad u_a(a)=u_a(a+c).
\end{align}
\end{proposition}

\begin{proof}
Let
\begin{align*}
    A(u,v):=\int_0^1 u'v'\,dx+\beta_0u(0)v(0)+\beta_1u(1)v(1),\qquad
    B_a(u,v):=\int_0^1 m_a uv\,dx .
\end{align*}
The eigenvalue equation is
\(A(u_a,v)=\lambda(a)B_a(u_a,v)\) for all \(v\in H^1(0,1)\), with
\(B_a(u_a,u_a)=1\). Since the eigenvalue branch is simple, the eigenfunction can be
chosen differentiably in \(a\) after this normalisation. Differentiating the weak
eigenvalue identity at \(v=u_a\), and using the symmetry of \(A\) and \(B_a\), gives
\begin{align*}
    2A(\partial_a u_a,u_a)
    &=\lambda'(a)B_a(u_a,u_a)+\lambda(a)\partial_aB_a(u_a,u_a)
    +2\lambda(a)B_a(\partial_a u_a,u_a).
\end{align*}
The terms involving \(\partial_a u_a\) cancel because
\(A(\partial_a u_a,u_a)=\lambda(a)B_a(\partial_a u_a,u_a)\), and hence
\begin{align*}
    \lambda'(a)=-\lambda(a)\partial_aB_a(u_a,u_a).
\end{align*}
It remains only to compute the derivative of the moving-weight form. In one dimension
all \(H^1\)-functions have continuous representatives. Therefore, for every continuous
\(w\),
\begin{align*}
    \frac{1}{h}\int_0^1\bigl(m_{a+h}-m_a\bigr)w^2\,dx
    \longrightarrow (\kappa+1)\bigl(w(a+c)^2-w(a)^2\bigr)
\end{align*}
as \(h\to0\). Applying this to \(w=u_a\), equivalently using
\(\partial_a m_a=(\kappa+1)(\delta_{a+c}-\delta_a)\) in the sense of distributions,
yields
\begin{align*}
    \partial_aB_a(u_a,u_a)=(\kappa+1)\bigl(u_a(a+c)^2-u_a(a)^2\bigl).
\end{align*}
This proves \eqref{eq:shape-derivative}. Since \(u_a>0\), the equality
\(u_a(a)^2=u_a(a+c)^2\) is equivalent to \(u_a(a)=u_a(a+c)\), and
\eqref{eq:critical-endpoint-values} follows.
\end{proof}

\begin{proposition}[Critical-point equation]\label{prop:critical-g}
Let \(a\in(0,1-c)\). Then the following conditions are equivalent:
\begin{align*}
    \lambda'(a)=0,\qquad u_a(a)=u_a(a+c),\qquad
    g(a,\beta_0,\beta_1,\lambda(a))=0.
\end{align*}
If \(\lambda(a)\ne\beta_0\beta_1\), this is equivalent to
\begin{align}\label{eq:a-star}
    a=\frac{1}{2\sqrt{\lambda(a)}}
    \tanh^{-1}\left(\frac{\sqrt{\lambda(a)}(\beta_1-\beta_0)}{\lambda(a)-\beta_0\beta_1}\right)
    +\frac{1-c}{2},
\end{align}
provided the argument of \(\tanh^{-1}\) lies in \((-1,1)\). If
\(\lambda(a)=\beta_0\beta_1\), then an interior critical point can occur only when
\(\beta_0=\beta_1\).
\end{proposition}

\begin{proof}
The equivalence between \(\lambda'(a)=0\) and \(u_a(a)=u_a(a+c)\) is
Proposition~\ref{prop:shape-derivative}. Since \(\lambda(a)\) is a simple zero of the
characteristic equation, differentiating
\(f(a,\beta_0,\beta_1,\kappa,\lambda(a))=0\) gives
\begin{align*}
    \partial_a f(a,\beta_0,\beta_1,\kappa,\lambda(a))
    +\partial_\lambda f(a,\beta_0,\beta_1,\kappa,\lambda(a))\lambda'(a)=0 .
\end{align*}
Thus \(\lambda'(a)=0\) is equivalent to \(\partial_a f=0\). By
Lemma~\ref{lem:spectral-positive-window}, \(\sin(\sqrt{\lambda(a)\kappa}c)>0\), and
\(\cosh(\sqrt{\lambda(a)}(2a+c-1))>0\). Hence \eqref{eq:df-a-factorised} shows that
\(\partial_a f=0\) is equivalent to \(g(a,\beta_0,\beta_1,\lambda(a))=0\). Solving this
last equation gives \eqref{eq:a-star} when \(\lambda(a)\ne\beta_0\beta_1\). If
\(\lambda(a)=\beta_0\beta_1\), then \(g=\sqrt{\lambda(a)}(\beta_0-\beta_1)\), so
\(g=0\) requires \(\beta_0=\beta_1\).
\end{proof}

\begin{theorem}[One-dimensional placement criterion]\label{thm:onedim}
Let \(c\in(0,1)\), \(\kappa>0\), \(\beta_0\ge0\), and \(\beta_1\ge0\), with at least
one of \(\beta_0,\beta_1\) positive. Let \(\lambda(a)\) be the positive principal
eigenvalue corresponding to \(m_a\), \(a\in[0,1-c]\), and define
\begin{align*}
    \mathcal M:=\operatorname*{argmin}_{a\in[0,1-c]}\lambda(a),\qquad
    \mathcal C:=\{a\in(0,1-c):g(a,\beta_0,\beta_1,\lambda(a))=0\}.
\end{align*}
Then \(\mathcal M\) is nonempty and compact, and
\begin{align}\label{eq:minimizer-characterization}
    \mathcal M\subseteq \{0,1-c\}\cup \mathcal C .
\end{align}
More precisely:
\begin{itemize}
    \item if \(\lambda'(a)>0\) for all \(a\in(0,1-c)\), then \(\mathcal M=\{0\}\);
    \item if \(\lambda'(a)<0\) for all \(a\in(0,1-c)\), then \(\mathcal M=\{1-c\}\);
    \item if \(\mathcal C=\{a_*\}\), \(\lambda'(a)<0\) on \((0,a_*)\), and
    \(\lambda'(a)>0\) on \((a_*,1-c)\), then \(\mathcal M=\{a_*\}\);
    \item if \(\mathcal C=\{a_*\}\), \(\lambda'(a)>0\) on \((0,a_*)\), and
    \(\lambda'(a)<0\) on \((a_*,1-c)\), then
    \(\mathcal M\subseteq\{0,1-c\}\), with the smaller endpoint value deciding the
    global minimiser.
\end{itemize}
At every interior candidate \(a_*\in\mathcal C\), the explicit relation
\eqref{eq:a-star} holds whenever \(\lambda(a_*)\ne\beta_0\beta_1\). {The derivative criterion and the equation \(g=0\) are used only for interior placements. Endpoint minimisers are decided by direct comparison of \(\lambda(0)\) and \(\lambda(1-c)\).}
\end{theorem}

\begin{proof}
Continuity of \(\lambda\) on the compact interval \([0,1-c]\) gives nonemptiness and
compactness of \(\mathcal M\). If a minimiser lies in \((0,1-c)\), Fermat's rule and
Lemma~\ref{lem:lambda-regularity} give \(\lambda'(a)=0\). By
Proposition~\ref{prop:critical-g}, this is equivalent to
\(g(a,\beta_0,\beta_1,\lambda(a))=0\), proving \eqref{eq:minimizer-characterization}.
The four subsequent alternatives are the corresponding elementary consequences of
the sign of the derivative on the intervals separated by critical points.
\end{proof}
\endgroup

\paragraph{{Remark on the equal-endpoint Robin case.}}
\begingroup\color{black}
{When the endpoint parameters are equal, \(\beta_0=\beta_1=\beta\), the criticality factor becomes
\[
    g(a,\beta,\beta,\lambda)
    =(\lambda-\beta^2)\tanh\left(2\sqrt{\lambda}\left(a-\frac{1-c}{2}\right)\right).
\]
Thus \(g=0\) is not equivalent, in general, to \(\lambda=\beta^2\): it also holds at the centred placement \(a=(1-c)/2\), independently of whether \(\lambda=\beta^2\). The reduction \(g=0\Rightarrow \lambda=\beta^2\) is valid only after excluding the centred critical point, or in the particular non-centred elimination used in the equal-Robin classification of \cite{Hintermueller2011}. Under asymmetric Robin parameters, the result above gives a rigorous branchwise characterisation: an interior minimiser must satisfy the coupled conditions \(f=0\) and \(g=0\), equivalently the endpoint-value condition \(u_a(a)=u_a(a+c)\). This avoids any unsupported spectral-window assumption and follows the actual branch \(\lambda=\lambda(a)\).}
\endgroup

\begingroup\color{black}
\paragraph{On the scope of the branchwise criterion.}
The branchwise criterion above is the form used in the remainder of the paper. Recovering a closed six-case statement in the asymmetric setting would require a separate global
monotonicity theorem for
\(G(a):=g(a,\beta_0,\beta_1,\lambda(a))\). The derivative of \(G\) contains the
additional term \(\partial_\lambda g\,\lambda'(a)\), and this term is not controlled by the
fixed-\(\lambda\) argument alone. At zeros of \(G\) the extra
term vanishes, which is useful for identifying the local crossing direction, but it does
not by itself give a global classification of all placements. We therefore keep the
rigorous endpoint/interior-candidate characterisation and treat the numerical section
as an illustration of this branchwise criterion.
\endgroup

\section{A Schr\"odinger-type extension with a fixed background potential}\label{sec:potential-extension}

\subsection{General variational and bang--bang results}\label{subsec:q-general}

The ecological interpretation of \eqref{eq:pb} suggests a natural extension in which the sign-changing weight $m$ is still the optimisable resource, but an additional spatially varying hostile field is prescribed in advance. More precisely, let $q\in L^\infty(\Omega)$ and consider
\begin{equation}\label{eq:pbq}
	\left\{
	\begin{aligned}
		\Delta \phi + (\lambda m(x)-q(x))\phi &= 0 && \text{in } \Omega,\\
		\partial_{\mathbf n}\phi + \beta\phi &= 0 && \text{on } \partial\Omega.
	\end{aligned}
	\right.
\end{equation}
When $q\equiv 0$, this reduces to \eqref{eq:pb}. We think of $q$ as a fixed background penalty: depending on the application, it may represent mortality, harvesting pressure, pollution, predation, or some other adverse environmental effect that is not itself subject to optimisation. In terms of the diffusive logistic model, the corresponding parabolic equation is
\begin{equation}\label{eq:logistic-q}
	\begin{cases}
		\displaystyle u_t = \Delta u + \bigl(\omega m(x)-q(x)\bigr)u-u^2 & \text{in } \Omega\times \mathbb R_+,\\[0.4ex]
		\partial_{\mathbf n}u + {\beta}u = 0 & \text{on } \partial\Omega\times \mathbb R_+,
	\end{cases}
\end{equation}
whose linearisation at $u\equiv 0$ is exactly \eqref{eq:pbq}. Thus the optimisation problem asks how a favourable region should be placed when the habitat is already subject to a non-removable hostile background. This fits naturally with the reaction--diffusion literature on heterogeneous environments, boundary loss, and harvesting or protection effects; see, for example, \cite{CantrellCosner2003,CantrellCosner2006,CuiLiMeiShi2017}. From the spectral point of view, \eqref{eq:pbq} is a Schr\"odinger-type problem: equivalently, one may rewrite it as $-\Delta \phi + q\phi = \lambda m\phi$ with Robin boundary condition. For related indefinite Robin eigenvalue problems on general domains, see also \cite{Daners2013}.

{We record a variational extension and a perturbative stability statement for the case with a fixed background potential. We do not claim a placement classification for general positive potential; instead, we keep the part of the theory that survives with essentially the same proofs as in the case $q\equiv 0$.} {In particular, all placement statements for the Schr\"odinger-type model below are limited to the one-dimensional constant-potential compactness result.} We restrict ourselves to the coercive situation
\begin{equation}\label{eq:coercive-assumption}
	q\ge 0\quad\text{a.e. in }\Omega,
	\qquad
	\beta\ge 0\quad\text{on }\partial\Omega,
\end{equation}
and assume that either $q>0$ on a subset of $\Omega$ of positive measure or $\beta>0$ on a subset of $\partial\Omega$ of positive surface measure. This avoids the degenerate pure Neumann case $(q,\beta)\equiv (0,0)$, where the sign of $\int_\Omega m$ again becomes decisive.

For $\phi\in H^1(\Omega)$ and $\lambda\in \mathbb R$, set
\[
Q^q_{\phi}(\lambda)
:=
\int_\Omega |\nabla \phi|^2
+\int_\Omega q\phi^2
+\int_{\partial\Omega} \beta\phi^2
-\lambda\int_\Omega m\phi^2,
\]
and
\[
\mu_q(\lambda)
:=
\inf_{\phi\in H^1(\Omega)}
\frac{Q^q_{\phi}(\lambda)}{\int_\Omega \phi^2}.
\]
Exactly as in Lemma~\ref{lem:smoller} and Lemma~\ref{lem:mu-prime}, one shows that $\lambda$ is a principal eigenvalue of \eqref{eq:pbq} if and only if $\mu_q(\lambda)=0$, and that if $\phi_0$ is a minimiser in the definition of $\mu_q(\lambda)$, then
\[
\mu_q'(\lambda)
=
-\frac{\int_\Omega m\phi_0^2}{\int_\Omega \phi_0^2}.
\]
Since these arguments are literally the same as before, we omit the proof.

\begin{theorem}\label{thm:q-principal}
	Assume \eqref{eq:coercive-assumption}, and assume further that either $q>0$ on a subset of $\Omega$ of positive measure or $\beta>0$ on a subset of $\partial\Omega$ of positive surface measure. Let $m\in L^\infty(\Omega)$ be sign-changing with
	\[
	-1\le m(x)\le \kappa\qquad\text{for a.e. }x\in \Omega.
	\]
	Then problem \eqref{eq:pbq} admits a unique positive principal eigenvalue, denoted by $\lambda_q(m)$. Moreover,
	\begin{equation}\label{eq:rayleigh-q}
		\lambda_q(m)
		=
		\inf_{\phi\in S(m)}
		\frac{\displaystyle
			\int_\Omega |\nabla \phi|^2
			+\int_\Omega q\phi^2
			+\int_{\partial\Omega} \beta\phi^2}
		{\displaystyle \int_\Omega m\phi^2},
		\qquad
		S(m):=\left\{\phi\in H^1(\Omega);\ \int_\Omega m\phi^2>0\right\}.
	\end{equation}
\end{theorem}

\begin{proof}
	For every fixed $\phi\in H^1(\Omega)$, the map $\lambda\mapsto Q^q_{\phi}(\lambda)$ is affine, hence concave. Therefore $\mu_q$ is concave as the infimum of affine functions. Since $m$ changes sign, there exist $\phi_1,\phi_2\in H^1(\Omega)$ such that
	\[
	\int_\Omega m\phi_1^2>0,
	\qquad
	\int_\Omega m\phi_2^2<0,
	\]
	and consequently $\mu_q(\lambda)\to -\infty$ both as $\lambda\to +\infty$ and as $\lambda\to -\infty$.
	
	Next consider
	\[
	\Lambda_1(q,\beta)
	:=
	\inf_{\phi\in H^1(\Omega)\setminus\{0\}}
	\frac{\displaystyle
		\int_\Omega |\nabla \phi|^2
		+\int_\Omega q\phi^2
		+\int_{\partial\Omega} \beta\phi^2}
	{\displaystyle \int_\Omega \phi^2}.
	\]
	By the compact embedding $H^1(\Omega)\hookrightarrow L^2(\Omega)$, the infimum is attained. Since $q\ge 0$ and $\beta\ge 0$, the numerator is nonnegative. If $\Lambda_1(q,\beta)=0$, then for a minimiser $\phi$ we must have
	\[
	\nabla\phi=0\quad\text{a.e. in }\Omega,
	\qquad
	q\phi=0\quad\text{a.e. in }\Omega,
	\qquad
	\beta\phi=0\quad\text{a.e. on }\partial\Omega.
	\]
	Hence $\phi$ is constant, and the additional assumption {that either $q>0$ on a subset of $\Omega$ of positive measure or $\beta>0$ on a subset of $\partial\Omega$ of positive surface measure} forces this constant to be zero, a contradiction. Therefore $\Lambda_1(q,\beta)>0$, and so
	\[
	\mu_q(0)\ge \Lambda_1(q,\beta)>0.
	\]
	By concavity, $\mu_q$ has exactly two zeros, one positive and one negative. The positive zero is the unique positive principal eigenvalue of \eqref{eq:pbq}, which we denote by $\lambda_q(m)$.
	
	Finally, since $\mu_q\bigl(\lambda_q(m)\bigr)=0$, the same argument as in Theorem~\ref{thm:principal} yields the Rayleigh quotient formula \eqref{eq:rayleigh-q}.
\end{proof}

We now consider the optimisation problem
\begin{align}\label{eq:min-problem-q}
	\inf_{m\in \mathcal M_{m_0,\kappa}} \lambda_q(m),
	\qquad
	\mathcal M_{m_0,\kappa}
	:=
	\left\{
	\begin{array}{l}
	m\in L^\infty(\Omega):\ -1\le m\le \kappa,\ |\Omega_m^+|>0,\\[0.2ex]
	\displaystyle \int_\Omega m\le -m_0|\Omega|
	\end{array}
	\right\}.
\end{align}
where $m_0\in (-\kappa,1)$.

\begin{theorem}\label{thm:q-bangbang}
	Under the assumptions of Theorem~\ref{thm:q-principal}, the optimisation problem \eqref{eq:min-problem-q} admits a minimiser. {Moreover, every minimiser is bang--bang: if $m^*\in\mathcal M_{m_0,\kappa}$ minimises \eqref{eq:min-problem-q}, then there exists a measurable set $E^*\subset\Omega$ such that}
	\begin{equation}\label{eq:q-bangbang}
		m^* = \kappa\chi_{E^*}-\chi_{\Omega\setminus E^*}
		\qquad\text{a.e. in }\Omega.
	\end{equation}
	{Every minimiser also satisfies the active volume constraint}
	\begin{equation}\label{eq:q-active}
		\int_\Omega m^* = -m_0|\Omega|.
	\end{equation}
	Equivalently, if
	\[
	c:=\frac{1-m_0}{1+\kappa},
	\]
	then the optimisation problem can be rewritten as
	\begin{equation}\label{eq:q-shape-problem}
		\inf_{E\subset \Omega,\ |E|=c|\Omega|} \lambda_q(E),
		\qquad
		\lambda_q(E):=\lambda_q\bigl(\kappa\chi_E-\chi_{\Omega\setminus E}\bigr).
	\end{equation}
\end{theorem}

\begin{proof}
	The existence part follows by the same direct-method argument as in \cite[Theorem~1.4]{Lou2006} and \cite[Theorem~1.2]{Hintermueller2011}. Indeed, compared with the case $q\equiv0$, the only modification is the additional fixed term $\int_\Omega q\phi^2$ in the numerator of the Rayleigh quotient, which is independent of the optimisation variable $m$ and does not affect the compactness or lower-semicontinuity arguments.

	{We prove the stronger bang--bang assertion for an arbitrary minimiser. Let $\widehat m$ be any minimiser, and let $\phi>0$ be an associated eigenfunction. Define
	\[
	A_q(\phi):=\int_\Omega |\nabla \phi|^2+\int_\Omega q\phi^2+\int_{\partial\Omega}\beta\phi^2,
	\qquad
	D_m(\phi):=\int_\Omega m\phi^2.
	\]
	Then $D_{\widehat m}(\phi)>0$ and $\lambda_q(\widehat m)=A_q(\phi)/D_{\widehat m}(\phi)$. For every admissible $m$ with $D_m(\phi)>0$,
	\[
	\lambda_q(\widehat m)\le \lambda_q(m)
	\le \frac{A_q(\phi)}{D_m(\phi)}.
	\]
	Thus $D_m(\phi)\le D_{\widehat m}(\phi)$, so $\widehat m$ maximises $m\mapsto D_m(\phi)$ among admissible competitors with positive denominator.

	The mass constraint is active. Indeed, if $\int_\Omega \widehat m<-m_0|\Omega|$, then one may increase $\widehat m$ by a sufficiently small positive amount on a subset of positive measure where $\widehat m<\kappa$. The new competitor is still admissible and still has positive denominator, while $D_m(\phi)$ strictly increases because $\phi>0$ a.e.; this contradicts maximality. Hence
	\[
	\int_\Omega\widehat m=-m_0|\Omega|.
	\]

	Applying the bathtub principle to $m\mapsto D_m(\phi)$ under the bounds $-1\le m\le\kappa$ and the active mass constraint gives a maximiser
	\[
	m^\sharp=\kappa\chi_{E^\sharp}-\chi_{\Omega\setminus E^\sharp},
	\qquad |E^\sharp|=c|\Omega|.
	\]
	Moreover $D_{m^\sharp}(\phi)\ge D_{\widehat m}(\phi)>0$. Therefore
	\[
	\lambda_q(\widehat m)
	\le\lambda_q(m^\sharp)
	\le \frac{A_q(\phi)}{D_{m^\sharp}(\phi)}
	\le \frac{A_q(\phi)}{D_{\widehat m}(\phi)}
	=\lambda_q(\widehat m).
	\]
	All inequalities are equalities. Hence $m^\sharp$ is also a minimiser and $\phi$ attains the Rayleigh quotient for $m^\sharp$, so $\phi$ is an eigenfunction for $m^\sharp$ with the same eigenvalue. The two weak equations are
	\[
	\Delta\phi+(\lambda_q(\widehat m)\widehat m-q)\phi=0,
	\qquad
	\Delta\phi+(\lambda_q(\widehat m)m^\sharp-q)\phi=0.
	\]
	Subtracting gives $(m^\sharp-\widehat m)\phi=0$ a.e. in $\Omega$. Since the principal eigenfunction is positive a.e. in $\Omega$, $\widehat m=m^\sharp$ a.e. Thus the arbitrary minimiser $\widehat m$ is bang--bang and satisfies the active mass condition. The identity \eqref{eq:q-active} then gives $|E^*|=c|\Omega|$, and \eqref{eq:q-shape-problem} follows.}
\end{proof}

\subsection{The one-dimensional interval family with constant background potential}\label{subsec:q0-interval}

{We now specialise \eqref{eq:pbq} to $\Omega=(0,1)$ and to a constant background potential $q\equiv q_0\ge 0$. We do not attempt a full analogue of Theorem~\ref{thm:onedim}; instead, we record a perturbative result which gives a rigorous bridge between the unperturbed problem and the Schr\"odinger-type model, and which can be illustrated numerically.}

Fix $c\in (0,1)$, let
\[
E_a:=(a,a+c),\qquad a\in [0,1-c],
\]
and define
\[
m_a:=\kappa\chi_{E_a}-\chi_{(0,1)\setminus E_a}.
\]
For every $a\in [0,1-c]$ and every $q_0\ge 0$, consider
\begin{equation}\label{eq:pbq-1d}
	\left\{
	\begin{aligned}
		\phi'' + (\lambda m_a(x)-q_0)\phi &= 0 && \text{in } (0,1),\\
		\phi'(0)-\beta_0\phi(0) &= 0,\\
		\phi'(1)+\beta_1\phi(1) &= 0,
	\end{aligned}
	\right.
\end{equation}
where {$\beta_0,\beta_1\ge 0$ and $\beta_0+\beta_1>0$}. By Theorem~\ref{thm:q-principal}, for every such pair $(a,q_0)$ the problem \eqref{eq:pbq-1d} admits a unique positive principal eigenvalue; we denote it by $\lambda_{q_0}(a)$.

The key observation is that the transfer--matrix formulation used in Section~\ref{sec:shooting} extends with essentially no change in strategy.

\begin{proposition}\label{prop:q0-transfer}
	For $m\in\{-1,\kappa\}$, set
	\[
	Q_m^{(q_0)}(\lambda):=
	\begin{pmatrix}
		0 & 1\\
		q_0-\lambda m & 0
	\end{pmatrix}.
	\]
	Then \eqref{eq:pbq-1d} is equivalent to the characteristic equation
	\begin{equation}\label{eq:q0-shooting}
		R_{q_0}(a,\lambda):=e_2M_{q_0}(a,\lambda)(1,\beta_0)^\top+\beta_1 e_1M_{q_0}(a,\lambda)(1,\beta_0)^\top=0,
	\end{equation}
	where
	\begin{equation}\label{eq:q0-transfer}
		M_{q_0}(a,\lambda):=e^{(1-a-c)Q_{-1}^{(q_0)}(\lambda)}e^{cQ_{\kappa}^{(q_0)}(\lambda)}e^{aQ_{-1}^{(q_0)}(\lambda)}.
	\end{equation}
	Moreover, $\lambda_{q_0}(a)$ is the smallest positive root of \eqref{eq:q0-shooting}.
\end{proposition}

\begin{proof}
	Let $w=(\phi,\phi')^\top$. On every subinterval where $m_a$ is constant, the equation in \eqref{eq:pbq-1d} becomes
	\[
	w'=Q_m^{(q_0)}(\lambda)w,
	\]
	and therefore the solution is propagated by matrix exponentials. Since \eqref{eq:pbq-1d} is homogeneous in $\phi$, we may normalise $\phi(0)=1$, so that the left Robin condition gives $w(0)=(1,\beta_0)^\top$. Propagating first across $(0,a)$, then $(a,a+c)$, and finally $(a+c,1)$ yields \eqref{eq:q0-transfer}. The right Robin condition is exactly \eqref{eq:q0-shooting}. The last assertion follows because the principal eigenvalue is the first eigenvalue of this regular one-dimensional problem.
\end{proof}

The next result is the basic perturbative statement needed for the Schr\"odinger extension.

\begin{proposition}\label{prop:q0-continuity}
	Fix $\overline q>0$. Then the map
	\[
	(a,q_0)\longmapsto \lambda_{q_0}(a)
	\]
	is continuous on $[0,1-c]\times [0,\overline q]$. In particular, for every $q_0\in [0,\overline q]$ the minimiser set
	\[
	\mathcal A(q_0):=\operatorname*{argmin}_{a\in [0,1-c]}\lambda_{q_0}(a)
	\]
	is nonempty and compact.
\end{proposition}

\begin{proof}
	By Proposition~\ref{prop:q0-transfer}, the characteristic function $R_{q_0}(a,\lambda)$ is real-analytic in $(a,q_0,\lambda)$. The same Lagrange-identity argument used in Lemma~\ref{lem:lambda-regularity}, with the equation $y''+(\lambda m_a-q_0)y=0$, gives
\[
 y(1)\,\partial_\lambda R_{q_0}(a,\lambda_{q_0}(a))=-\int_0^1 m_a y^2\,dx.
\]
For the positive principal eigenvalue the integral is positive and $y(1)>0$, so the principal zero is simple. The implicit function theorem therefore yields local continuous dependence of $\lambda_{q_0}(a)$ on $(a,q_0)$, and hence continuity on the whole compact rectangle $[0,1-c]\times [0,\overline q]$.
	
	The nonemptiness and compactness of $\mathcal A(q_0)$ then follow from the continuity of $a\mapsto \lambda_{q_0}(a)$ on the compact interval $[0,1-c]$.
\end{proof}

\begin{theorem}\label{thm:q0-limit-minimisers}
	Let $q_j\to 0^+$, and for each $j$ choose some $a_j\in \mathcal A(q_j)$. Then every accumulation point of $(a_j)$ belongs to $\mathcal A(0)$. Equivalently, every family of minimisers for the constant-potential problem converges, along subsequences, to minimisers of the unperturbed problem.
\end{theorem}

\begin{proof}
	Since $a_j\in [0,1-c]$, the sequence $(a_j)$ has an accumulation point; passing to a subsequence if necessary, we may assume that $a_j\to a_*\in [0,1-c]$. Let $a\in [0,1-c]$ be arbitrary. Since each $a_j$ minimises $\lambda_{q_j}$, we have
	\[
	\lambda_{q_j}(a_j)\le \lambda_{q_j}(a).
	\]
	Now let $j\to\infty$. By Proposition~\ref{prop:q0-continuity},
	\[
	\lambda_0(a_*)=\lim_{j\to\infty}\lambda_{q_j}(a_j)
	\le \lim_{j\to\infty}\lambda_{q_j}(a)=\lambda_0(a).
	\]
	Since $a\in [0,1-c]$ was arbitrary, this shows that $a_*\in \mathcal A(0)$.
\end{proof}

{Combining Theorem~\ref{thm:q0-limit-minimisers} with the branchwise characterisation in Theorem~\ref{thm:onedim}, we obtain the following consequence.}

\begin{corollary}\label{cor:q0-from-unperturbed}
\begingroup\color{black}
Let \(q_j\to0^+\), and for each \(j\) choose some \(a_j\in\mathcal A(q_j)\). If
\(\mathcal A(0)\) denotes the minimiser set of the unperturbed problem, then
\begin{align}\label{eq:q0-distance-general}
    \operatorname{dist}(a_j,\mathcal A(0))\to0 .
\end{align}
In particular, if \(\mathcal A(0)=\{a_0\}\), then \(a_j\to a_0\).
\endgroup
\end{corollary}

\begin{proof}
\begingroup\color{black}
If \eqref{eq:q0-distance-general} failed, then there would be an \(\varepsilon>0\) and a
subsequence, still denoted \(a_j\), such that
\(\operatorname{dist}(a_j,\mathcal A(0))\ge\varepsilon\). By compactness of \([0,1-c]\),
we may pass to a further subsequence with \(a_j\to a_*\). Theorem~\ref{thm:q0-limit-minimisers}
then gives \(a_*\in\mathcal A(0)\), contradicting the distance bound.
\endgroup
\end{proof}

{Corollary~\ref{cor:q0-from-unperturbed} is a compactness statement describing how minimisers for the constant-potential problem approach the unperturbed minimiser set as \(q_0\downarrow0\). It gives a clean target for the numerical part: one should observe the minimisers of \(a\mapsto\lambda_{q_0}(a)\) drifting towards the unperturbed minimiser set characterised by Theorem~\ref{thm:onedim}.}

\section{{Numerical illustrations by a transfer--matrix shooting method}}
\label{sec:shooting}

{This section collects the numerical illustrations accompanying the two one-dimensional parts of the paper. We first illustrate the unperturbed transfer-matrix criteria, and then illustrate Corollary~\ref{cor:q0-from-unperturbed} for the constant-potential perturbation.} In both cases the computation is based on the same transfer--matrix shooting principle. {The adaptive-shooting archive is in \cite{dataset}.}

\subsection{{Numerical illustration of the unperturbed one-dimensional criteria}}\label{subsec:illustrate-thm-onedim}

{We now recompute the unperturbed one-dimensional examples using the branchwise criterion of Theorem~\ref{thm:onedim}. The purpose is no longer to establish or check a closed six-case formula. Instead, the computation illustrates the endpoint-value diagnostic supplied by the shape derivative.}

For given $(c,\kappa,\beta_0,\beta_1)$ and a placement $a\in[0,1-c]$ of the favourable subinterval $(a,a+c)$, the equation
\begin{align}\label{eq:pb-1d-restated}
	\begin{cases}
		u''+\lambda\,m(x)\,u=0,\\
		u'(0)-\beta_0 u(0)=0,\\
		u'(1)+\beta_1 u(1)=0,
	\end{cases}
\end{align}
with 
\begin{align}
	m(x)=\begin{cases} \kappa,&\text{ if }x\in(a,a+c),\\ -1,&\text{ otherwise,}\end{cases}
\end{align}
admits explicit fundamental solutions on each subinterval. Writing the state vector as $w=(u,u')^\top$, we have
\begin{align}
	w'(x)=\begin{pmatrix}0&1\\-\lambda m(x)&0 \end{pmatrix}w(x)=:Q_{m(x)}(\lambda)w(x).
\end{align}
The continuity of $u$ and $u'$ at $x=a$ and $x=a+c$ yields
\begin{align}\label{eq:transfer}
	w(1)=e^{(1-a-c)Q_{-1}(\lambda)} e^{cQ_\kappa(\lambda)} e^{aQ_{-1}(\lambda)} w(0).
\end{align}
After the normalisation $u(0)=1$, the left Robin condition gives $w(0)=(1,\beta_0)^\top$, and the right Robin condition gives the scalar residual equation
\begin{align}\label{eq:shooting}
	R(a,\lambda):=e_2 M(a,\lambda)\,(1,\beta_0)^\top+\beta_1\,e_1 M(a,\lambda)\,(1,\beta_0)^\top = 0,
\end{align}
where $M(a,\lambda):=e^{(1-a-c)Q_{-1}(\lambda)} e^{cQ_\kappa(\lambda)} e^{aQ_{-1}(\lambda)}$. The principal eigenvalue is computed as the smallest positive root of \eqref{eq:shooting}.

{The root search is adaptive. In all computations reported below we use the lower bound $\lambda_{\min}=10^{-12}$, the initial upper bound $\lambda_{\max}^{(0)}=1$, enlargement factor $2$, and $N_\lambda=500$ uniform scan subintervals at each enlargement step. After the first sign-changing interval of $R(a,\lambda)$ is found, it is refined by bisection until the interval length is below $10^{-10}$. This avoids relying on any invalid fixed window such as $(0,\pi^2/(4c^2\kappa))$. The plotted profiles use the same placement grid as the tables.}

{For each sampled placement $a_j$, we also compute the shape-derivative diagnostic}
\begin{align}\label{eq:numerical-shape-derivative-diagnostic}
    D(a_j):=u_{a_j}(a_j)^2-u_{a_j}(a_j+c)^2.
\end{align}
{By Proposition~\ref{prop:shape-derivative}, the sign of $D(a_j)$ is the sign of $\lambda'(a_j)$, up to the positive factor $\lambda(a_j)(\kappa+1)$ under the normalisation used there. Thus a sign change from negative to positive indicates an interior minimum candidate, while a sign change from positive to negative indicates an interior maximum candidate.}

{We take $c=0.5$, $\kappa=2$, and sample $N_a=80$ equal subintervals of $[0,1-c]$. Table~\ref{tab:unperturbed-illustration} records four canonical behaviours of the branchwise criterion. The labels ``left'', ``right'', ``interior'', and ``two endpoint'' are numerical descriptions of the computed minimiser set, not cases of a closed explicit classification theorem. The values of $D$ are shown at representative points in the left and right halves of the interval.}

\begin{table}[!htbp]
	\centering
	\caption{{Adaptive unperturbed computation illustrating the branchwise criterion of Theorem~\ref{thm:onedim}.}}
	\label{tab:unperturbed-illustration}
	\scriptsize
	\setlength{\tabcolsep}{2pt}
	\begin{tabular}{ccccccc}
		\toprule
		case & $\beta_0$ & $\beta_1$ & $\operatorname{argmin}\lambda$ & $\lambda_{\min}$ & $D_L$ & $D_R$ \\
		\midrule
		unique left endpoint & $0.2$ & $5.0$ & $\{0\}$ & $1.597391$ & $0.535472$ & $0.949119$ \\
		unique right endpoint & $5.0$ & $0.2$ & $\{1-c\}$ & $1.597391$ & $-3.028538$ & $-9.474211$ \\
		unique interior point & $10.0$ & $10.0$ & $\{(1-c)/2\}$ & $5.288192$ & $-3.971210$ & $13.790605$ \\
		two endpoint minimisers & $0.2$ & $0.2$ & $\{0,1-c\}$ & $0.631873$ & $0.066311$ & $-0.084110$ \\
		\bottomrule
	\end{tabular}
\end{table}

\begin{figure}[!htbp]
	\centering
	\begin{minipage}[t]{0.48\textwidth}\centering
		\includegraphics[width=\linewidth]{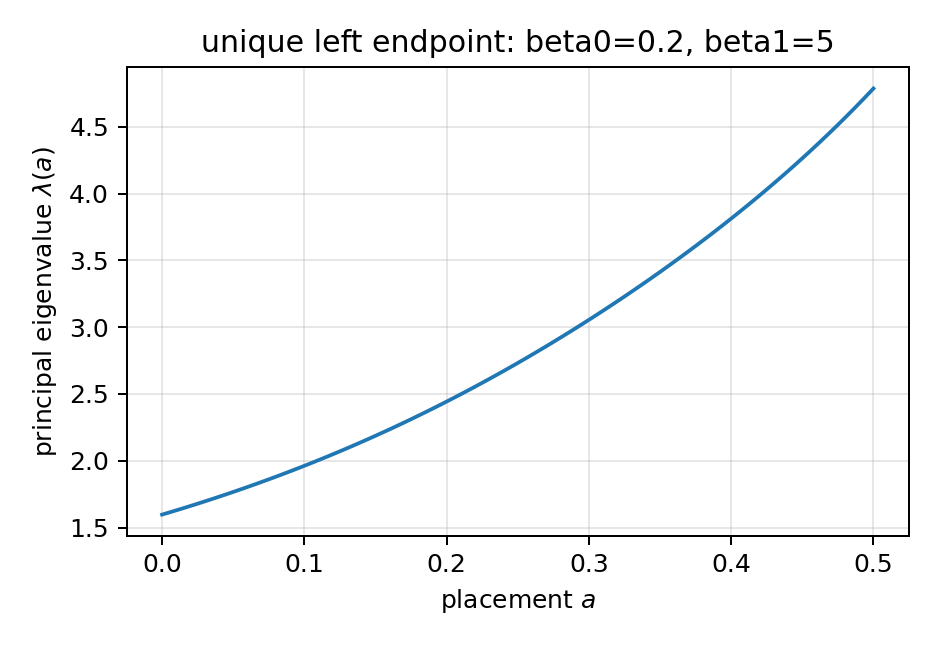}
		\\[-2pt]{\small unique left endpoint}
	\end{minipage}\hfill
	\begin{minipage}[t]{0.48\textwidth}\centering
		\includegraphics[width=\linewidth]{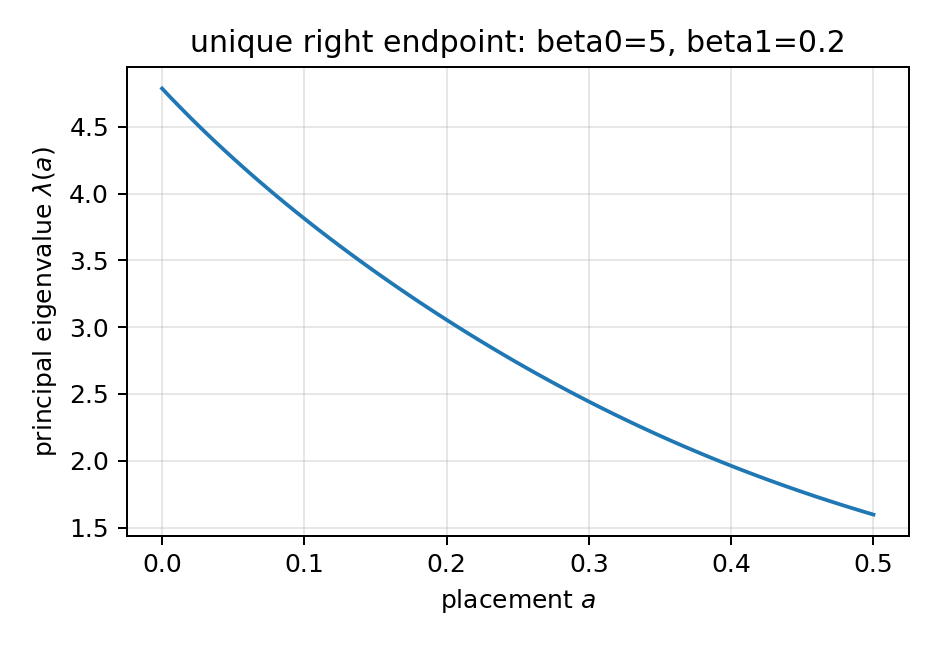}
		\\[-2pt]{\small unique right endpoint}
	\end{minipage}
	
	\medskip
	\begin{minipage}[t]{0.48\textwidth}\centering
		\includegraphics[width=\linewidth]{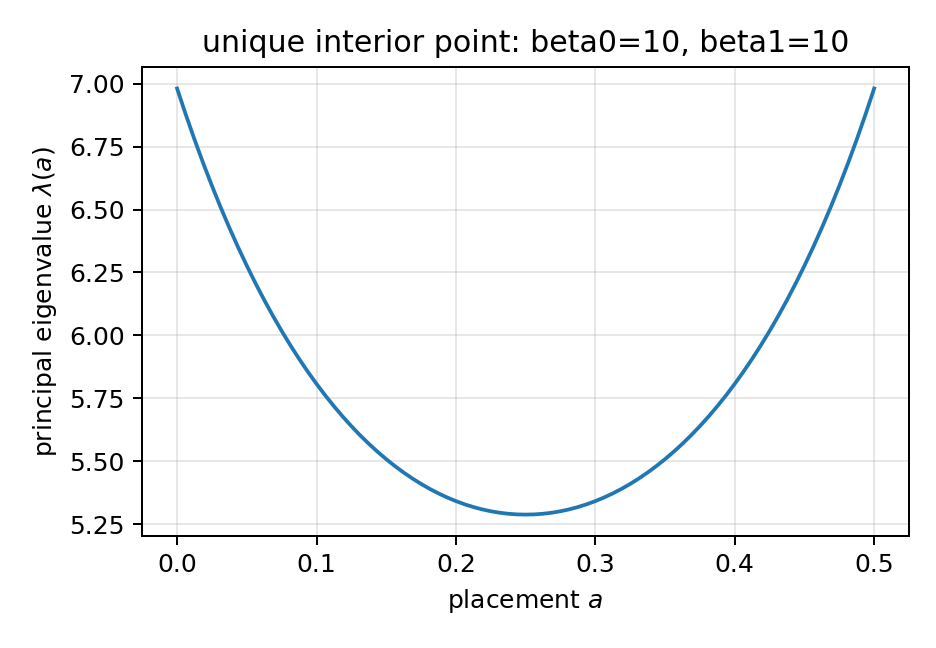}
		\\[-2pt]{\small unique interior point}
	\end{minipage}\hfill
	\begin{minipage}[t]{0.48\textwidth}\centering
		\includegraphics[width=\linewidth]{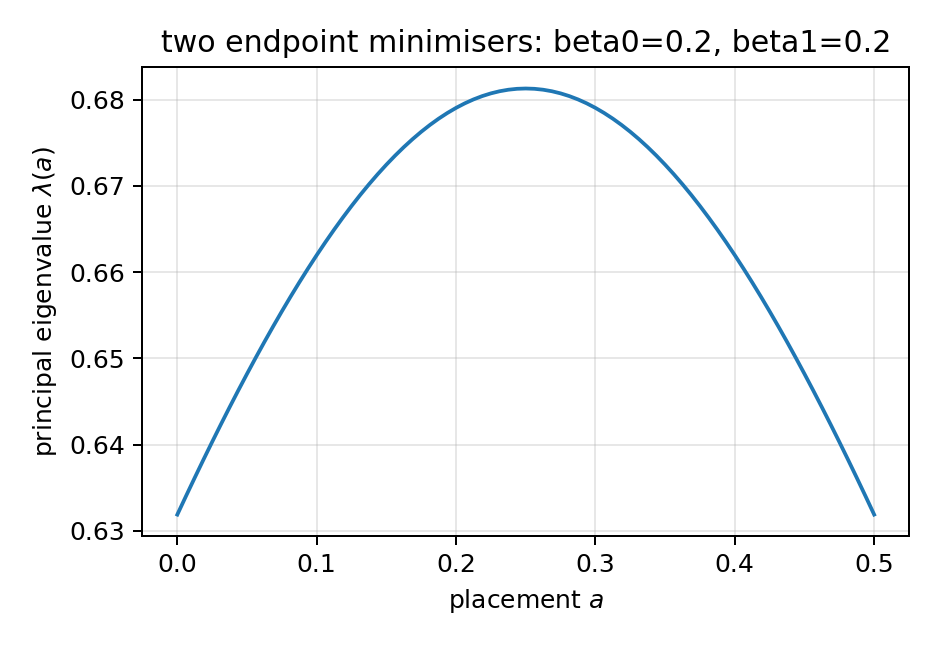}
		\\[-2pt]{\small two endpoint minimisers}
	\end{minipage}
	\caption{{Profiles of $a\mapsto\lambda(a)$ for the four adaptive unperturbed computations in Table~\ref{tab:unperturbed-illustration}.}}
	\label{fig:unperturbed-profiles}
\end{figure}

\subsection{{Numerical illustration of Corollary~\ref{cor:q0-from-unperturbed}}}
\label{subsec:q0-illustration}

{We now illustrate Corollary~\ref{cor:q0-from-unperturbed} by the same transfer--matrix strategy as above.} For fixed $q_0\ge 0$ and $a\in[0,1-c]$, Proposition~\ref{prop:q0-transfer} reduces the eigenvalue problem \eqref{eq:pbq-1d} to the scalar equation
\[
R_{q_0}(a,\lambda)=0,
\]
and the principal eigenvalue $\lambda_{q_0}(a)$ is computed as the smallest positive root. Numerically, this is done by {the same adaptive scan-and-bisection procedure, with the same values of $\lambda_{\min}$, $\lambda_{\max}^{(0)}$, enlargement factor, scan size, and bisection tolerance as in Subsection~\ref{subsec:illustrate-thm-onedim}}. The only modification with respect to Section~\ref{sec:shooting} is that the matrices $Q_m$ are replaced by
\[
Q_m^{(q_0)}(\lambda)=\begin{pmatrix}0&1\\ q_0-\lambda m&0\end{pmatrix}.
\]

For the placement variable, we work on a uniform grid $\{a_j\}_{j=0}^{N_a}\subset[0,1-c]$ and define the discrete minimiser set
\[
\mathcal A_{N_a}(q_0):=\arg\min\{\lambda_{q_0}(a_j):0\le j\le N_a\}.
\]
To compare with Corollary~\ref{cor:q0-from-unperturbed}, we also record the deterministic selector
\[
a^{\sharp}(q_0):=\min \mathcal A_{N_a}(q_0),
\]
that is, the leftmost point in the discrete argmin set. If the unperturbed minimiser is unique, then $a^{\sharp}(q_0)$ should approach it as $q_0\downarrow 0$; in the two-endpoint case, the distance from $a^{\sharp}(q_0)$ to $\{0,1-c\}$ should tend to $0$.

\begin{figure}[!htbp]
	\centering
	\begin{minipage}[t]{0.48\textwidth}\centering
		\includegraphics[width=\linewidth]{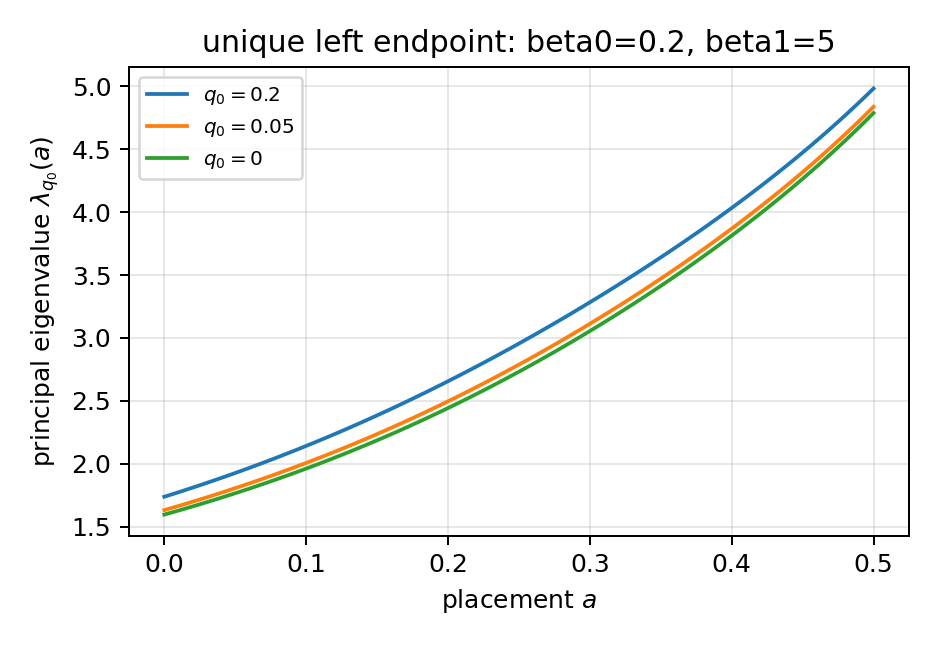}
		\\[-2pt]{\small $(\beta_0,\beta_1)=(0.2,5.0)$: left endpoint.}
	\end{minipage}\hfill
	\begin{minipage}[t]{0.48\textwidth}\centering
		\includegraphics[width=\linewidth]{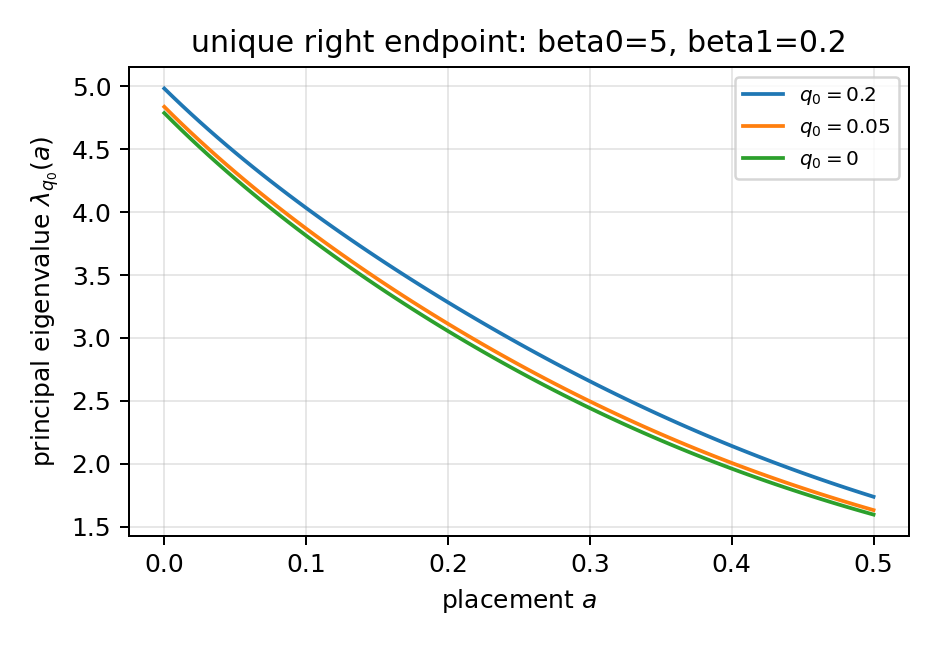}
		\\[-2pt]{\small $(\beta_0,\beta_1)=(5.0,0.2)$: right endpoint.}
	\end{minipage}
	
	\medskip
	\begin{minipage}[t]{0.48\textwidth}\centering
		\includegraphics[width=\linewidth]{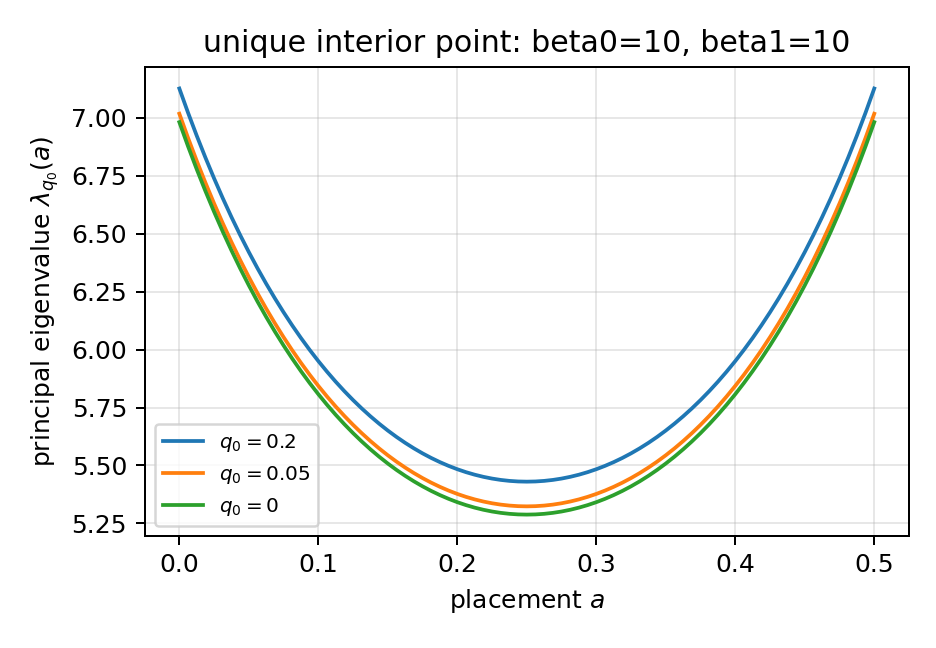}
		\\[-2pt]{\small $(\beta_0,\beta_1)=(10.0,10.0)$: unique interior point.}
	\end{minipage}\hfill
	\begin{minipage}[t]{0.48\textwidth}\centering
		\includegraphics[width=\linewidth]{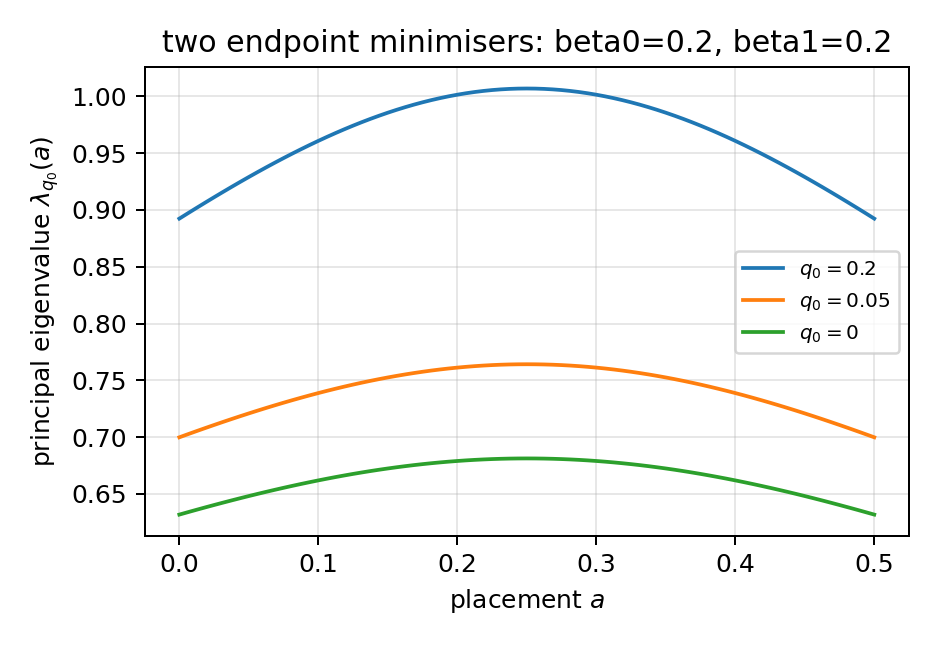}
		\\[-2pt]{\small $(\beta_0,\beta_1)=(0.2,0.2)$: two endpoint minimisers.}
	\end{minipage}
	\caption{{Profiles of $a\mapsto \lambda_{q_0}(a)$ for the four canonical tests in Table~\ref{tab:q0-illustration}. In each panel, the three curves correspond to $q_0=0.2$, $q_0=0.05$, and $q_0=0$.}}
	\label{fig:q0-profiles}
\end{figure}

\paragraph{Choice of parameters.}
{We use the same four unperturbed canonical choices as in Table~\ref{tab:unperturbed-illustration}:}
\[
    c=0.5,\qquad \kappa=2.0,\qquad
    (\beta_0,\beta_1)\in\{(0.2,5.0),(5.0,0.2),(10.0,10.0),(0.2,0.2)\}.
\]
{The placement grid has $N_a=80$ subintervals. The background-potential values used for the table are}
\[
q_0\in\{0.2,0.1,0.05,0.02,0.01,0\},
\]
{and Figure~\ref{fig:q0-profiles} shows the three representative curves $q_0=0.2$, $q_0=0.05$, and $q_0=0$.}

\begin{table}[!htbp]
	\centering
	\caption{{Adaptive numerical illustration of Corollary~\ref{cor:q0-from-unperturbed} for four canonical parameter choices. In each row, the displayed discrete argmin set is unchanged for every sampled value $q_0\in\{0.2,0.1,0.05,0.02,0.01,0\}$.}}
	\label{tab:q0-illustration}
	\scriptsize
	\setlength{\tabcolsep}{2pt}
	\begin{tabular}{ccccc}
		\toprule
		case & $(\beta_0,\beta_1)$ & $\mathcal A_{N_a}(0)$ & $\mathcal A_{N_a}(q_0)$ & max. distance \\
		\midrule
		unique left endpoint & $(0.2,5.0)$ & $\{0\}$ & always $\{0\}$ & $0$ \\
		unique right endpoint & $(5.0,0.2)$ & $\{1-c\}=\{0.5\}$ & always $\{0.5\}$ & $0$ \\
		unique interior point & $(10.0,10.0)$ & $\{(1-c)/2\}=\{0.25\}$ & always $\{0.25\}$ & $0$ \\
		two endpoint minimisers & $(0.2,0.2)$ & $\{0,1-c\}=\{0,0.5\}$ & always $\{0,0.5\}$ & $0$ \\
		\bottomrule
	\end{tabular}
\end{table}

{Table~\ref{tab:q0-illustration} illustrates Corollary~\ref{cor:q0-from-unperturbed}. For these representative examples, the discrete minimiser set is already stable over the sampled range of $q_0$; in less symmetric or near-transition regimes, the theorem only predicts convergence to the unperturbed minimiser set as $q_0\downarrow0$.}

{Together, Table~\ref{tab:q0-illustration} and Figure~\ref{fig:q0-profiles} provide a compact numerical illustration of the perturbative compactness statement.} {The purpose of Table~\ref{tab:q0-illustration} is only to provide a sanity check of stability in representative regimes; it is not intended as a numerical demonstration of visible drift of minimisers.} {The code producing these plots is included in the computational archive and uses the same adaptive root-search principle as the unperturbed computation.}

\section{Conclusion}

{We have studied the minimisation of the positive principal eigenvalue for an indefinite-weight problem with asymmetric Robin parameters. After recalling the variational framework and the standard bang--bang reduction, we analysed the one-dimensional placement problem for the optimal favourable interval. The one-dimensional result gives a rigorous branchwise characterisation: a global minimiser is either an endpoint placement or an interior point satisfying the coupled transfer-matrix criticality condition \(g(a,\beta_0,\beta_1,\lambda(a))=0\), equivalently \(u_a(a)=u_a(a+c)\). The argument follows the actual branch \(\lambda=\lambda(a)\) and avoids imposing an unsupported fixed spectral window.}

{We also introduced a Schr\"odinger-type extension in which the optimisable indefinite weight is coupled with a fixed nonnegative background potential. In the coercive setting, we proved that the principal-eigenvalue theory and the bang--bang structure persist. For the one-dimensional constant-potential family, we established a perturbative compactness result showing that minimisers for small positive potential converge, along subsequences, to minimisers of the unperturbed problem. The numerical examples illustrate this behaviour for representative configurations using the same adaptive root-search protocol.} {We do not claim a placement classification for general positive background potentials.}

{These results are relevant from two complementary viewpoints. On the modelling side, they describe how favourable habitat may be placed when the two ends of the environment interact differently with the exterior and when an additional hostile background field is present. On the analytical side, they furnish a concrete bridge between indefinite-weight optimisation problems arising in reaction--diffusion equations and Schr\"odinger-type spectral optimisation. The transfer--matrix formulation also yields a reproducible computational framework that can be adapted to related one-dimensional optimisation problems.}

Several directions remain open. The most immediate one is to move beyond the perturbative constant-potential regime and seek a sharper description of the optimal placement for positive background potential. It would also be natural to investigate whether the present asymmetric Robin analysis can be combined with the change-of-variables techniques available for shells and other special geometries. More broadly, the higher-dimensional shape problem remains largely open, and the present one-dimensional results provide a useful benchmark for future analytical and numerical work.

\section*{Statements and Declarations}

\subsection*{Funding}
This work was co-funded by the University of Ostrava, Grant No.~SGS05/PŘF/2026. Yifan Zhang is also co-funded by the European Union under the REFRESH -- Research Excellence For REgion Sustainability and High-tech Industries project, project number CZ.10.03.01/00/22\_003/0000048, via the Operational Programme Just Transition.

\subsection*{Competing interests}
The authors have no competing interests to declare that are relevant to the content of this article.

\subsection*{Ethics approval and consent to participate}
Not applicable

\subsection*{Consent for publication}
Not applicable

\subsection*{Data availability}
{The computational archive, including the adaptive-shooting scripts, machine-readable tables, and figures used in Section~\ref{sec:shooting}, is available at \href{https://doi.org/10.5281/zenodo.19219633}{10.5281/zenodo.19219633}.}

\subsection*{Materials availability}
Not applicable

\subsection*{Code availability}
{The code used to generate the datasets is included in the computational archive available at \href{https://doi.org/10.5281/zenodo.19219633}{10.5281/zenodo.19219633}.}

\subsection*{Author contribution}
All authors contributed to the mathematical development of the paper, to the interpretation of the results, and to the preparation of the manuscript. Yifan Zhang carried out the numerical experiments, {implemented the computational experiments}, and prepared the associated numerical data. All authors reviewed the manuscript and approved the final version.


\begin{thebibliography}{99}

\bibitem{Afrouzi1999}
Afrouzi, G. A.; Brown, K. J.: On principal eigenvalues for boundary value problems with indefinite weight and {R}obin boundary conditions. \emph{Proceedings of the American Mathematical Society} \textbf{127}(1), 125--130 (1999). DOI: \href{https://doi.org/10.1090/s0002-9939-99-04561-x}{10.1090/s0002-9939-99-04561-x}.

\bibitem{Afrouzi_2002}
Afrouzi, G. A.: Boundedness and monotonicity of principal eigenvalues for boundary value problems with indefinite weight functions. \emph{International Journal of Mathematics and Mathematical Sciences} \textbf{30}(1), 25--29 (2002). DOI: \href{https://doi.org/10.1155/s0161171202007780}{10.1155/s0161171202007780}.

\bibitem{Bocher1914}
B{\^o}cher, Maxime: The smallest characteristic numbers in a certain exceptional case. \emph{Bulletin of the American Mathematical Society} \textbf{21}(1), 6--9 (1914). DOI: \href{https://doi.org/10.1090/s0002-9904-1914-02560-1}{10.1090/s0002-9904-1914-02560-1}.

\bibitem{BerestyckiNirenbergVaradhan1994}
Berestycki, Henri; Nirenberg, Louis; Varadhan, S. R. S.: The principal eigenvalue and maximum principle for second-order elliptic operators in general domains. \emph{Communications on Pure and Applied Mathematics} \textbf{47}(1), 47--92 (1994). DOI: \href{https://doi.org/10.1002/cpa.3160470105}{10.1002/cpa.3160470105}.

\bibitem{Brown1980}
K. J. Brown; S. S. Lin: On the existence of positive eigenfunctions for an eigenvalue problem with indefinite weight function. \emph{Journal of Mathematical Analysis and Applications} \textbf{75}(1), 112--120 (1980). DOI: \href{https://doi.org/10.1016/0022-247x(80)90309-1}{10.1016/0022-247x(80)90309-1}.

\bibitem{Cantrell_1989}
Cantrell, Robert Stephen; Cosner, Chris: Diffusive logistic equations with indefinite weights: population models in disrupted environments. \emph{Proceedings of the Royal Society of Edinburgh: Section A Mathematics} \textbf{112}(3--4), 293--318 (1989). DOI: \href{https://doi.org/10.1017/s030821050001876x}{10.1017/s030821050001876x}.

\bibitem{CantrellCosner2003}
Cantrell, Robert Stephen; Cosner, Chris: \emph{Spatial Ecology via Reaction-Diffusion Equations}. John Wiley \& Sons, Ltd., Chichester (2003). DOI: \href{https://doi.org/10.1002/0470871296}{10.1002/0470871296}.

\bibitem{CantrellCosner2006}
Cantrell, Robert Stephen; Cosner, Chris: On the effects of nonlinear boundary conditions in diffusive logistic equations on bounded domains. \emph{Journal of Differential Equations} \textbf{231}(2), 768--804 (2006). DOI: \href{https://doi.org/10.1016/j.jde.2006.08.018}{10.1016/j.jde.2006.08.018}.

\bibitem{Clarte2021}
Clart{\'e}, Thibaut T.; Schaeffer, Nathana{\"e}l; Labrosse, St{\'e}phane; Vidal, J{\'e}r{\'e}mie: The effects of a {R}obin boundary condition on thermal convection in a rotating spherical shell. \emph{Journal of Fluid Mechanics} \textbf{918} (2021). DOI: \href{https://doi.org/10.1017/jfm.2021.356}{10.1017/jfm.2021.356}.

\bibitem{CuiLiMeiShi2017}
Cui, Renhao; Li, Haomiao; Mei, Linfeng; Shi, Junping: Effect of harvesting quota and protection zone in a reaction-diffusion model arising from fishery management. \emph{Discrete and Continuous Dynamical Systems - Series B} \textbf{22}(3), 791--807 (2017). DOI: \href{https://doi.org/10.3934/dcdsb.2017039}{10.3934/dcdsb.2017039}.

\bibitem{Daners2000}
Daners, Daniel: Robin boundary value problems on arbitrary domains. \emph{Transactions of the American Mathematical Society} \textbf{352}(9), 4207--4236 (2000). DOI: \href{https://doi.org/10.1090/S0002-9947-00-02444-2}{10.1090/S0002-9947-00-02444-2}.

\bibitem{Daners2013}
Daners, Daniel: Principal eigenvalues for generalised indefinite Robin problems. \emph{Potential Analysis} \textbf{38}(4), 1047--1069 (2013). DOI: \href{https://doi.org/10.1007/s11118-012-9306-9}{10.1007/s11118-012-9306-9}.

\bibitem{Hess1991}
Hess, Peter: \emph{Periodic-Parabolic Boundary Value Problems and Positivity}. Longman Scientific \& Technical, Harlow (1991).

\bibitem{Hintermueller2011}
Hinterm{\"u}ller, M.; Kao, C.-Y.; Laurain, A.: Principal eigenvalue minimization for an elliptic problem with indefinite weight and {R}obin boundary conditions. \emph{Applied Mathematics and Optimization} \textbf{65}(1), 111--146 ({2012}). DOI: \href{https://doi.org/10.1007/s00245-011-9153-x}{10.1007/s00245-011-9153-x}.

\bibitem{Kao2008}
Chiu-Yen Kao; Yuan Lou; Eiji Yanagida: Principal eigenvalue for an elliptic problem with indefinite weight on cylindrical domains. \emph{Mathematical Biosciences and Engineering} \textbf{5}(2), 315--335 (2008). DOI: \href{https://doi.org/10.3934/mbe.2008.5.315}{10.3934/mbe.2008.5.315}.

\bibitem{Krein1955}
Kre{\u{\i}}n, M. G.: On certain problems on the maximum and minimum of characteristic values and on the {Lyapunov} zones of stability. In: \emph{Eleven Papers on Topology, Function Theory and Differential Equations}, pp.~163--187. American Mathematical Society, Providence (1955). DOI: \href{https://doi.org/10.1090/trans2/001/08}{10.1090/trans2/001/08}.

\bibitem{Lamboley2016}
Lamboley, Jimmy; Laurain, Antoine; Nadin, Gr{\'e}goire; Privat, Yannick: Properties of optimizers of the principal eigenvalue with indefinite weight and {R}obin conditions. \emph{Calculus of Variations and Partial Differential Equations} \textbf{55}(6) (2016). DOI: \href{https://doi.org/10.1007/s00526-016-1084-6}{10.1007/s00526-016-1084-6}.

\bibitem{Lou2006}
Lou, Yuan; Yanagida, Eiji: Minimization of the principal eigenvalue for an elliptic boundary value problem with indefinite weight, and applications to population dynamics. \emph{Japan Journal of Industrial and Applied Mathematics} \textbf{23}(3), 275--292 (2006). DOI: \href{https://doi.org/10.1007/bf03167595}{10.1007/bf03167595}.

\bibitem{MunozCastaneda2015}
Mu\~noz-Casta\~neda, J. M.; Mateos Guilarte, J.: $\delta$-$\delta'$ generalized {R}obin boundary conditions and quantum vacuum fluctuations. \emph{Physical Review D} \textbf{91}(2), 025028 (2015). DOI: \href{https://doi.org/10.1103/physrevd.91.025028}{10.1103/physrevd.91.025028}.

\bibitem{Murray2002}
Murray, J. D.: \emph{Mathematical Biology. {I}. An Introduction}. Springer, New York, 3 ed. (2002). DOI: \href{https://doi.org/10.1007/b98868}{10.1007/b98868}.

\bibitem{OkuboLevin2001}
Okubo, Akira; Levin, Simon A.: \emph{Diffusion and Ecological Problems: Modern Perspectives}. Springer, New York, 2 ed. (2001). DOI: \href{https://doi.org/10.1007/978-1-4757-4978-6}{10.1007/978-1-4757-4978-6}.

\bibitem{Schneider_2025}
Schneider, Baruch; Schneiderov{\'a}, Diana; Zhang, Yifan: Optimization of {R}obin {L}aplacian eigenvalue with indefinite weight in spherical shell. \emph{Mathematical Methods in the Applied Sciences} \textbf{48}(6), 6586--6591 (2025). DOI: \href{https://doi.org/10.1002/mma.10697}{10.1002/mma.10697}.

\bibitem{Senn1982}
Senn, Stefan; Hess, Peter: On positive solutions of a linear elliptic eigenvalue problem with neumann boundary conditions. \emph{Mathematische Annalen} \textbf{258}(4), 459--470 (1982). DOI: \href{https://doi.org/10.1007/bf01453979}{10.1007/bf01453979}.

\bibitem{Skellam_1951}
Skellam, J. G.: Random dispersal in theoretical populations. \emph{Biometrika} \textbf{38}(1/2), 196 (1951). DOI: \href{https://doi.org/10.2307/2332328}{10.2307/2332328}.

\bibitem{Smoller_1994}
Smoller, Joel: \emph{Shock Waves and Reaction--Diffusion Equations}. Springer, New York, 2 ed. (1994). DOI: \href{https://doi.org/10.1007/978-1-4612-0873-0}{10.1007/978-1-4612-0873-0}.

\bibitem{dataset}
{Zhang, Yifan: Code and data for one-dimensional optimisation of indefinite-weight principal eigenvalues with asymmetric Robin parameters and a Schr{\"o}dinger-type perturbation, adaptive-shooting archive.} Zenodo (2026). {DOI: \href{https://doi.org/10.5281/zenodo.19219633}{10.5281/zenodo.19219633}.}

\end{thebibliography}
\end{document}